\documentclass[12pt]{amsart}

\usepackage[bbgreekl]{mathbbol}
\usepackage{mathrsfs}
\usepackage{graphicx}
\usepackage{amsmath}
\usepackage{amsfonts}
\usepackage{indentfirst}
\usepackage{amsthm}
\usepackage{amssymb}
\usepackage{algorithm, algorithmic}
\usepackage{xcolor}
\usepackage{geometry}
\usepackage{appendix}
\usepackage{bm}
\usepackage{url}
\usepackage{setspace}
\usepackage{subfigure}
\usepackage{verbatim}
\usepackage[normalem]{ulem}
\usepackage{multirow}
\usepackage{lineno}
\usepackage{float}
\usepackage{array}
\usepackage{diagbox}

\newtheorem{theorem}{Theorem}[section]
\newtheorem{lemma}{Lemma}[section]
\newtheorem{remark}{Remark}[section]

\renewcommand{\theequation}{\arabic{section}.\arabic{equation}}

\renewcommand{\thelemma}{\arabic{section}.\arabic{lemma}}

\renewcommand{\thetable}{\arabic{section}.\arabic{table}}

\newcommand\tbbint{{-\mkern -16mu\int}}

\newcommand\dbbint{{-\mkern -19mu\int}}

\newcommand\bbint{
{\mathchoice{\dbbint}{\tbbint}{\tbbint}{\tbbint}}
}

\date{\today}

\begin{document}

\title{QM/MM Methods for Crystalline Defects. Part 3: 
	Machine-Learned Interatomic Potentials
}

\author{Huajie Chen}
\address{School of Mathematical Sciences, Beijing Normal University, Beijing 100875, China} 
\email{chen.huajie@bnu.edu.cn}
\author{Christoph Ortner}
\address{University of British Columbia, 1984 Mathematics Road, Vancouver, BC, Canada}
\email{ortner@math.ubc.ca}
\author{Yangshuai Wang}
\address{Corresponding author. Institute of Natural Sciences and School of Mathematical Sciences, Shanghai Jiao Tong University, Shanghai 200240, China}
\email{yswang2016@sjtu.edu.cn}


\begin{abstract}
We develop and analyze a framework for consistent QM/MM (quantum/classic) hybrid models of crystalline defects, which admits general atomistic interactions including traditional off-the-shell interatomic potentials as well as state of art ``machine-learned interatomic potentials''.
We (i) establish an {\it a priori} error estimate for the QM/MM approximations in terms of matching conditions between the MM and QM models, and (ii) demonstrate how to use these matching conditions to construct practical machine learned MM potentials specifically for QM/MM simulations.
%
%
\end{abstract}

\maketitle


\renewcommand\arraystretch{1.5}

\newcommand\R{\mathbb{R}}
\newcommand\Rg{\mathcal{R}}
\newcommand\N{\mathbb{N}}
\newcommand\Z{\mathbb{Z}}
\newcommand\C{\mathbb{C}}
\newcommand\E{\mathcal{E}}
\newcommand\n{\mathfrak{n}}
\def\asP{\textnormal{\bf (P)}}
\def\asD{\textnormal{\bf (D)}}
\def\asRC{\textnormal{\bf (RC)}}
\def\asSE{\textnormal{\bf (S)}}
\def\asSER{\textnormal{\bf (S.R)}}
\def\assERL{\textnormal{\bf (RL)}}
\def\asSEL{\textnormal{\bf (S.L)}}
\def\asSEH{\textnormal{\bf (S.H)}}
\def\asSEP{\textnormal{\bf (S.P)}}
\def\asSEI{\textnormal{\bf (S.I)}}
\def\asSEPS{\textnormal{\bf (S.PS)}}
\def\asLS{\textnormal{\bf (LS)}}
\def\asS{\textnormal{\bf (S)}}
\def\pt{\textnormal{\bf (PT)}}
\def\VMM{V^{\rm MM}}
\def\VQM{V^{\rm QM}}
\def\Vta{V^{\rm Taylor}}
\def\Vhom{V^{\rm h}}
\def\Ih{I_1^{\rm h}}
\def\Wb{\mathbf{W}}
\def\Vsh{V^{\rm SHIPs}}
\def\E{\mathcal{E}}
\def\EH{\E^{\rm H}}
\def\uh{u_{\rm h}}
\def\vh{v_{\rm h}}
\def\p0{\pmb{0}}
\def\Egfc{\E^{\rm GFC}}
\def\EHT{\widetilde{\E}^{\rm H}}
\def\Vb{V^{\rm BUF}_{\#}}
\def\Usx{\mathscr{U}^{\rm H}}
\def\Us{\mathscr{U}}
\def\uH{\bar{u}^{\rm H}}
\def\UsH{{\mathscr{U}}^{1,2}}
\newcommand\UsHd{\Us^{-1,2}}
\def\THu{T^{\rm H}\bar{u}}
\def\lgamma{\ell^2_{\gamma}}
\def\efit{\varepsilon^{\rm E}}
\def\ffit{\varepsilon^{\rm F}}
\def\vfit{\varepsilon^{\rm V}}
\def\L{\Lambda}
\def\R{\mathbb{R}}
\def\Rl{\mathcal{R}}
\def\Rg{\mathcal{R}^{*}}
\def\Adm{{\rm Adm}}
\newcommand{\<}{\langle}
\renewcommand{\>}{\rangle}
\def\LQM{\Lambda^{\rm QM}}
\def\LMM{\Lambda^{\rm MM}}
\def\LFF{\Lambda^{\rm FF}}
\def\Lbuf{\Lambda^{\rm BUF}}
\def\Lhom{\L^{\rm h}}
\def\LhomS{\L^{\rm h}_*}
\def\OQM{\Omega^{\rm QM}}
\def\OMM{\Omega^{\rm MM}}
\def\OFF{\Omega^{\rm FF}}
\def\Obuf{\Omega^{\rm BUF}}
\def\RQM{R_{\rm QM}}
\def\RMM{R_{\rm MM}}
\def\RFF{R_{\rm FF}}
\def\Rbuf{R_{\rm BUF}}
\def\rcut{R_{\rm c}}
\def\Rcore{R_{\rm DEF}}
\def\FMM{\F^{\rm MM}}
\def\FQM{\F^{\rm QM}}
\def\F{\mathcal{F}}
\def\N{\mathbb{N}}
\def\FH{\F^{\rm H}}
\def\Fa{\widetilde{\F}}
\def\Ea{\widetilde{\E}}
\def\Adm{{\rm Adm}}
\def\Admu{\mathscr{A}}
\def\wf{\mathfrak{w}}
\def\Hw{\mathscr{L}}
\def\fc{f_{\rm c}}
\def\Rc{R_{\rm c}}
\def\dx{\,{\rm d}x}
\def\dt{\,{\rm d}t}
\def\Wcb{W_{\rm cb}}
\def\A{\mathbf{A}}
\def\n{\mathfrak{n}}


\section{Introduction}
\label{sec:introduction}
\setcounter{equation}{0}

Quantum mechanics and molecular mechanics (QM/MM) coupling methods have been widely used in multiscale modeling and simulations of large systems in materials science and biology \cite{bernstein09,csanyi04,gao02,kermode08,sherwood08,zhang12}. 
The region of primary interest (e.g., a defect core) is described by a QM model, which is embedded in an ambient environment (e.g., bulk crystal) described by an MM model. 
In principle, QM/MM coupling methods promise (near-)QM accuracy at (near-)MM computational cost for large-scale atomistic simulations. To achieve this in practise the two key approximation parameters that must be controlled are the size of the QM region and the choice of MM model. The purpose of this paper is to devise a framework for constructing {\em minimal} MM models matching the QM model in a way that results in explicit convergence rates.


QM/MM methods for material simulations originate in the need for large-scale simulations due to long-ranged elastic or electric fields that are strongly coupled to material defects. For example, the long-range fields give rise to an interaction between far-away defects, and classical continuum theories are typically sufficient to describe this. On the other hand, the resultant effect of those long-range fields on the motion of a defect (e.g. diffusion of point defects towards a grain boundary, or dislocation dynamics) strongly depends on the core structure and details of the chemistry of bond breaking.
These general ideas only motivate our own work, hence we refer to other sources for deeper introductions and further references on the applications of QM/MM schemes \cite{bernstein09, csanyi04, kermode08, khare2007coupled, reuter2000frontier}.

Our own interest is in the question how to design QM/MM schemes in a way that allows rigorous error control. To that end, we focus on a simpler but instructive setting of a single defect embedded in a homogeneous host crystal, where a rigorous numerical analysis approach is feasible. The atomistic equilibration problem in this context is a well-defined variational problem, and properties of equilibrium configurations are precisely characterised~\cite{Ehrlacher16,chen18,chen19}. With this information in hand it is possible to specify conditions on MM models and QM/MM coupling schemes under which the schemes converges as the QM radius increases~\cite{bernstein09, CMAME, chen17, li2018pexsi, wang2020posteriori}. For example in \cite{chen17} it was shown that the rate of convergence in terms of the QM region radius is determined by clearly specified matching conditions on the Taylor expansion coefficients in the QM and MM interaction laws around a reference state. This observation is the starting point for the present work.

The key ingredient in the construction of a QM/MM scheme is the choice of MM model (given a high-fidelity QM reference model to be used for the QM region). Traditionally, empirical or even analytical interatomic potentials have been used in this context, since they more than adequately describe the elastic far-field, and are computationally inexpensive. The crux is that they generally do not match the QM interaction at all, which leads to the difficult problem of constructing coupling schemes that alleviate that mismatch. 

An entirely different approach, developed outside of the QM/MM context and in fact intended to replace QM/MM models altogether, are machine-learned interatomic potentials (MLIPs) \cite{behler07, zhang12, Braams09, bachmayr19, Drautz19, Shapeev16}. Their defining feature is that they employ universal approximators as their parameterisation and aim to approximate a reference QM potential energy landscape to within arbitrary accuracy on a large training set. We will leverage this technology in order to construct MM schemes that satisfy precise matching conditions identified in \cite{chen17} and guarantee a QM/MM scheme that is consistent with the reference pure QM model. 

To that end, we will first generalize the error analysis of \cite{chen17} to obtain error estimates in terms of a quantity that incorporates the matching conditions and, most importantly, can be easily incorporated into a loss function for training an MLIP. We will then train various ACE MLIPs~\cite{Drautz19,bachmayr19} using different variations of our general idea. The most promising candidates are those potentials where we replace fully atomistic matching conditions with a suitable continuum limit, while retaining a fully controllable error.

\subsubsection*{Outline}
In Section \ref{sec:equilibration}, we review a rigorous framework for modelling the geometry equilibration of crystalline defects. 
In Section \ref{sec:qmmm}, we construct a QM/MM coupling method for crystalline defects, and give {\it a priori} general error estimates for both energy-mixing and force-mixing schemes, with errors quantified in terms of certain matching conditions between the MM and QM models.
In Section \ref{sec:constructions}, we discuss how to apply our theory to construct an MLIP from the given QM model, and present the NRL tight binding model and ACE interatomic potential as examples that we will use in our tests.
In Section \ref{sec:numerics}, we perform the numerical experiments for some real materials with point defect and dislocation. 
Finally, some concluding remarks are given.
For simplicity of the presentations, the details of the proofs are given in appendices.

\subsubsection*{Notation}
For the sake of brevity, we will denote $A\backslash\{a\}$ by $A\backslash a$, and $\{b-a~\vert ~b\in A\}$ by $A-a$.
We will use the symbol $\langle\cdot,\cdot\rangle$ to denote an abstract duality pair between a Banach space and its dual space. The symbol $|\cdot|$ normally denotes the Euclidean or Frobenius norm, while $\|\cdot\|$ denotes an operator norm.
For a functional $E \in C^2(X)$, the first and second variations are denoted by $\<\delta E(u), v\>$ and $\<\delta^2 E(u) v, w\>$ for $u,v,w\in X$, respectively.
For $j\in\N$, ${\bm{g}}\in (\R^d)^A$, and $V \in C^j\big((\R^d)^A\big)$, we define the notation
\begin{eqnarray}\label{eq:dV}
	V_{,{\bm \rho}}\big({\bm g}\big) :=
	\frac{\partial^j V\big({\bm g}\big)}
	{\partial {\bm g}_{\rho_1}\cdots\partial{\bm g}_{\rho_j}}
	\qquad{\rm for}\quad{\bm \rho}=(\rho_1, \ldots, \rho_j)\in A^{j}.
	\end{eqnarray}
The symbol $C$ denotes generic positive constant that may change from one line of an estimate to the next. When estimating rates of decay or convergence, $C$ will always remain independent of the system size, the lattice configuration and the test functions. The dependence of $C$ will be clear from the context or stated explicitly.


\section{Equilibration of crystalline defects}
\label{sec:equilibration}
\setcounter{equation}{0}

A rigorous framework for modelling the geometric equilibrium of crystalline defects has been developed in \cite{Ehrlacher16} for finite-range interaction models (interatomic potentials) and extended in \cite{chen19} to a range of ab initio interaction models.
These two works formulate the equilibration of crystal defects as a variational problem in a discrete energy space and establish qualitatively sharp far-field decay estimates for the corresponding equilibrium configuration. 
We will review the framework as important background for our own work.

\subsection{Displacement space}
\label{sec:reference}

Let $d\in\{2,3\}$ be the (effective) dimension of the system.
We consider a single defect embedded in an infinite homogeneous crystalline bulk. 
Both point defects and straight line dislocations will be studied in this paper.

A homogeneous crystal reference configuration is given by the Bravais lattice
$\Lhom=A\Z^d$, for some non-singular matrix $A \in \mathbb{R}^{d \times d}$. 
The reference configuration for a system with defect is a set $\L \subset \R^d$ satisfying
\begin{flushleft}\label{as:asRC}
\asRC	\quad
$\exists ~\Rcore>0$, such that	$\L\backslash B_{\Rcore} = \Lhom\backslash B_{\Rcore}$ and $\L \cap B_{\Rcore}$ is finite.
\end{flushleft}

The deformed configuration of the infinite lattice $\L$ is a map $y:\L\rightarrow\R^d$.
We can decompose the configuration $y$ into
\begin{eqnarray}
\label{config_y_u}
y(\ell) = y_0(\ell) + u(\ell) = x_0(\ell) + u_0(\ell) + u(\ell)  \qquad\forall~\ell\in\Lambda,
\end{eqnarray}
where $x_0:\L\rightarrow\R^d$ is a linear map a crystalline reference configuration, $u_0: \L\rightarrow\R^d$ is a {\it far-field predictor} solving a continuum linearised elasticity (CLE) equation \cite{Ehrlacher16} enforcing the presence of the defect of interest and $u: \L\rightarrow\R^d$ is a {\it corrector}. For point defects we simply take $u_0(\ell)=0~\forall~\ell\in\L$. The derivation of $u_0$ for straight dislocations is reviewed in \ref{sec:appendixU0}. 


We collect the settings for point defects and dislocations in the following conditions \asP~and \asD~respectively.

\begin{flushleft}
	\asP \quad  \label{as:asP}
	$d \in \{2, 3 \}$;
	$x_0(\ell) = \ell~\forall~\ell\in\L$;
	$u_0(\ell)=0~\forall~\ell\in\L$.
\end{flushleft}
\begin{flushleft}
	\asD \quad   \label{as:asD}
	$\L=A\Z^2$;
	$x_0(\ell) $ given by \eqref{eq:P} and
	$u_0$ by \eqref{predictor-u_0-dislocation}.
\end{flushleft}

The set of possible atomic configurations is
\begin{align*}
\Adm_{0}(\L) &:= \bigcup_{\mathfrak{m}>0} \Adm_{\mathfrak{m}}(\L)
\qquad \text{with} \\
\Adm_{\mathfrak{m}}(\L) &:= \left\{ y:\L \rightarrow \R^{d}, ~
|y(\ell)-y(m)| > \mathfrak{m} |\ell-m|
\quad\forall~  \ell, m \in \L \right\},
\end{align*}
where the parameter $\mathfrak{m}>0$ prevents the accumulation of atoms.

For $\ell\in\L$ and $\rho\in\L-\ell$, we define the finite difference
$D_\rho u(\ell) := u(\ell+\rho) - u(\ell)$. For a subset $\Rl \subset \Lambda-\ell$, we
define $D_{\Rl} u(\ell) := (D_\rho u(\ell))_{\rho\in\Rl}$, and we consider $Du(\ell) := D_{\Lambda-\ell} u(\ell)$ to be a finite-difference stencil with infinite range.
For a stencil $Du(\ell)$, we define the stencil norms
\begin{align}\label{eq: nn norm}
\big|Du(\ell)\big|_{\mathcal{N}} := \bigg( \sum_{\rho\in \mathcal{N}(\ell) - \ell} \big|D_\rho u(\ell)\big|^2 \bigg)^{1/2}  
\quad{\rm and}\quad
\|Du\|_{\ell^2_{\mathcal{N}}} := \bigg( \sum_{\ell \in \L} |Du(\ell)|_{\mathcal{N}}^2 \bigg)^{1/2},
\end{align}
where $\mathcal{N}(\ell)$ is the set containing nearest neighbours of site $\ell$,
%
\begin{multline}
\label{def1:Nl:r}
\qquad
\mathcal{N}(\ell) := \Big\{ \, \tilde{A}\big(\ell\pm e_i\big), ~i=1,\cdots,d \, \Big| \, \tilde{A}\in\R^{d\times d} \text{ satisfies }
\\[1ex] 
\tilde{A}\Z^d = A\Z^d \text{ and } \|\tilde{A}\|_1 =\min_{B\in\R^{d\times d},~B\Z^d = A\Z^d}\|B\|_1 \, \Big\}. \,\,
\qquad
\end{multline}
We can then define the corresponding functional space of finite-energy displacements
\begin{align}\label{space:UsH}
\UsH(\L) := \big\{u:\L\rightarrow\mathbb{R}^{d} ~\big\lvert~ \|Du\|_{\ell^2_{\mathcal{N}}(\L)}<\infty \big\}
\end{align}
with the associated semi-norm $\|Du\|_{\ell^2_{\mathcal{N}}}$.
We also define the following subspace of compact displacements
\begin{align}\label{space:Uc}
\Us^{\rm c}(\L) := \big\{u:\L\rightarrow\mathbb{R}^{d} ~\big\lvert~ \exists~R >0~{\rm s.t.}~u = {\rm const}~{\rm in}~\L\setminus B_{R}\big\}.
\end{align}
Then the associated class of admissible displacements is given by
\begin{eqnarray*}
	\Admu(\L):= \big\{ u\in\UsH(\L) ~:~ y_0 +u\in\Adm_0(\L) \big\}.
\end{eqnarray*}


\subsection{The QM site potential}
\label{sec:siteE}
The site potential represents the local energy contributed from each atomic site.
It specifies the physical model used in the simulations. Interatomic potentials are usually {\em defined} in terms of a site potential, incorporating the modelling assumption that interaction are local. QM models are defined in terms of total energies, but a site energy can still be constructed in some cases. Here, we review the construction and analysis of \cite{chen16, finnis03}.

Let $\L$ be the reference configuration satisfying \asRC~with $\Lhom$ the corresponding homogeneous lattice.
Denote $\LhomS:= \Lhom \setminus 0 = \Lhom-0$.
We consider the site potential to be a collection of mappings $V_{\ell}:(\R^d)^{\L-\ell}\rightarrow\R$, which represent the energy distributed to each atomic site.
We make the following assumptions on the regularity and locality of the site potentials, which has been justified for some basic quantum mechanic models \cite{chen18,chen16,chen19tb, co2020}.
We refer to \cite[\S 2.3 and \S 4]{chen19} for discussions of more general site potentials.

\begin{itemize}
	\label{as:SE:pr}
	\item[\assERL]
	{\it Regularity and locality:}
	For all $\ell \in \L$, $V_{\ell}\big(Du(\ell)\big)$ possesses partial derivatives up to $\mathfrak{n}$-th order with $\n\geq 3$. For $j=1,\ldots,\n$, there exist constants $C_j$ and $\eta_j$ such that
	\begin{eqnarray}
	\label{eq:Vloc}
	\big|V_{\ell,{\bm \rho}}\big(Du(\ell)\big)\big|  \leq
	C_j \exp\Big(-\eta_j\sum^j_{l=1}|{\bm \rho}_l|\Big)
	\end{eqnarray} 
	for all $\ell \in \L$ and ${\bm \rho} \in (\L - \ell)^{j}$. 
\end{itemize}

If the reference configuration $\L$ is a homogeneous lattice, $\L = \Lhom$, then the site potential does not depend on site $\ell\in\Lhom$ due to the translation invariance. In this case, 
we will denote the site potential by $\Vhom:(\R^d)^{\LhomS}\rightarrow\R$.

Although the site potentials are defined on infinite stencils $(\R^d)^{\L-\ell}$, the setting also applies to finite systems or to finite range interactions. It is only necessary to assume in this case that the potential $V_{\ell}(\pmb{g})$ does not depend on the reference sites $\pmb{g}_{\rho}$ outside the interaction range.
In particular, we will denote by $V_{\ell}^{\Omega}$ the site potential of a finite system with the reference configuration lying in $\L\cap\Omega$.

\subsection{Geometry equilibration}
\label{sec:geometry}

Let $\L$ satisfy \asRC~and the site potential satisfy the assumptions \asSE. 
Then we have from \cite[Theorem 2.1]{chen19} that the energy-difference functional
\begin{eqnarray}
\label{energy-difference}
\E(u) := \sum_{\ell\in\Lambda}\Big(V_{\ell}\big(Du_0(\ell)+Du(\ell)\big)-V_{\ell}\big(Du_0(\ell)\big)\Big),
\end{eqnarray}
after an elementary renormalisation,
is well-defined on the admissible displacements set $\Admu(\L)$ and is $(\n-1)$-times continuously differentiable with respect to the $\|D\cdot\|_{\ell^2_{\mathcal{N}}}$ norm, for both point defects and dislocations.  

We can now rigorously formulate the variational problem for the equilibrium state as
\begin{equation}\label{eq:variational-problem}
\bar{u} \in \arg\min \big\{ \E(u), u \in \Admu \big\},
\end{equation}
where ``$\arg\min$'' is understood as the set of local minimisers.
The minimizer $\bar{u}$ satisfies the following first and second order optimality conditions
\begin{eqnarray}\label{eq:optimality}
\big\< \delta\E(\bar{u}) , v\big\> = 0 , \qquad \big\< \delta^2\E(\bar{u}) v , v\big\> \geq 0, \qquad\forall~v\in\UsH(\L) .
\end{eqnarray}
For the purpose of an approximation error analysis we will need a stronger stability condition, which could be proven but must be formulated as an {\em assumption}:
\begin{flushleft}\label{as:LS}
	\asS \quad {\it Strong stability:}
	$\exists~\bar{c}>0$ s.t. $
		\big\< \delta^2\E(\bar{u}) v , v\big\> \geq \bar{c} \|Dv\|^2_{\ell^2_{\mathcal{N}}} \qquad\forall~v\in\UsH(\L).
$
\end{flushleft}

The next result gives the decay estimates for the equilibrium state for point defects and dislocations \cite[Theorem 3.2 and 3.7]{chen19}: Suppose that either \asP~or \asD~is satisfied, if $\bar{u}\in\Admu(\L)$ is a strongly stable solution to \eqref{eq:variational-problem} satisfying \asS, then there exists $C > 0$ such that
\begin{eqnarray} 
\label{eq:ubar-decay}
\big|D\bar{u}(\ell)\big|_{\mathcal{N}} \leq C \big(1+|\ell|\big)^{-d}\log^{t}(2+|\ell|),
\end{eqnarray}
where $t=0$ for \asP~and $t=1$ for \asD.

Instead of energy minimization, we may alternatively consider the force equilibrium formulation (that corresponds to the first part of \eqref{eq:optimality}):
\begin{eqnarray}
\label{eq:problem-force}
{\rm Find}~ \bar{u} \in \Admu(\L), ~~ {\rm s.t.} \quad
\F_{\ell}(\bar{u}) = 0, \qquad \forall~\ell\in\Lambda ,
\end{eqnarray}
where $\F_{\ell}:= -\nabla_{\ell} \E(u)$ represents the force on the atomic site $\ell$. 
Note that any of the minimizers of \eqref{eq:variational-problem} also solve \eqref{eq:problem-force}. 
For $u^{\rm h}\in\Admu(\Lhom)$ on a homogeneous lattice $\Lhom$, the force on each atomic site $\ell$ satisfies $\F_{\ell}(u^{\rm h})=\F^{\rm h}\big(u^{\rm h}(\cdot-\ell)\big)$ with some homogeneous force $\F^{\rm h}$ that does not depend on site $\ell$. 
%


\section{QM/MM coupling methods with MLIPs}
\label{sec:qmmm}
\setcounter{equation}{0}

QM/MM coupling schemes can be categorized into energy- and force-mixing schemes: energy-based methods construct a hybrid total energy functional and seek a local minimizer of that functional, while force-based methods solve the force balance equation with QM and MM contributions. We will review both kinds of schemes in more detail, and then derive the main theoretical results of this paper.
Towards that end we first need to discuss the domain partitioning and provide a framework to quantify how well an MM model matches the reference QM model.

\begin{figure}[!htb]
        \centering
        \includegraphics[height=6.0cm]{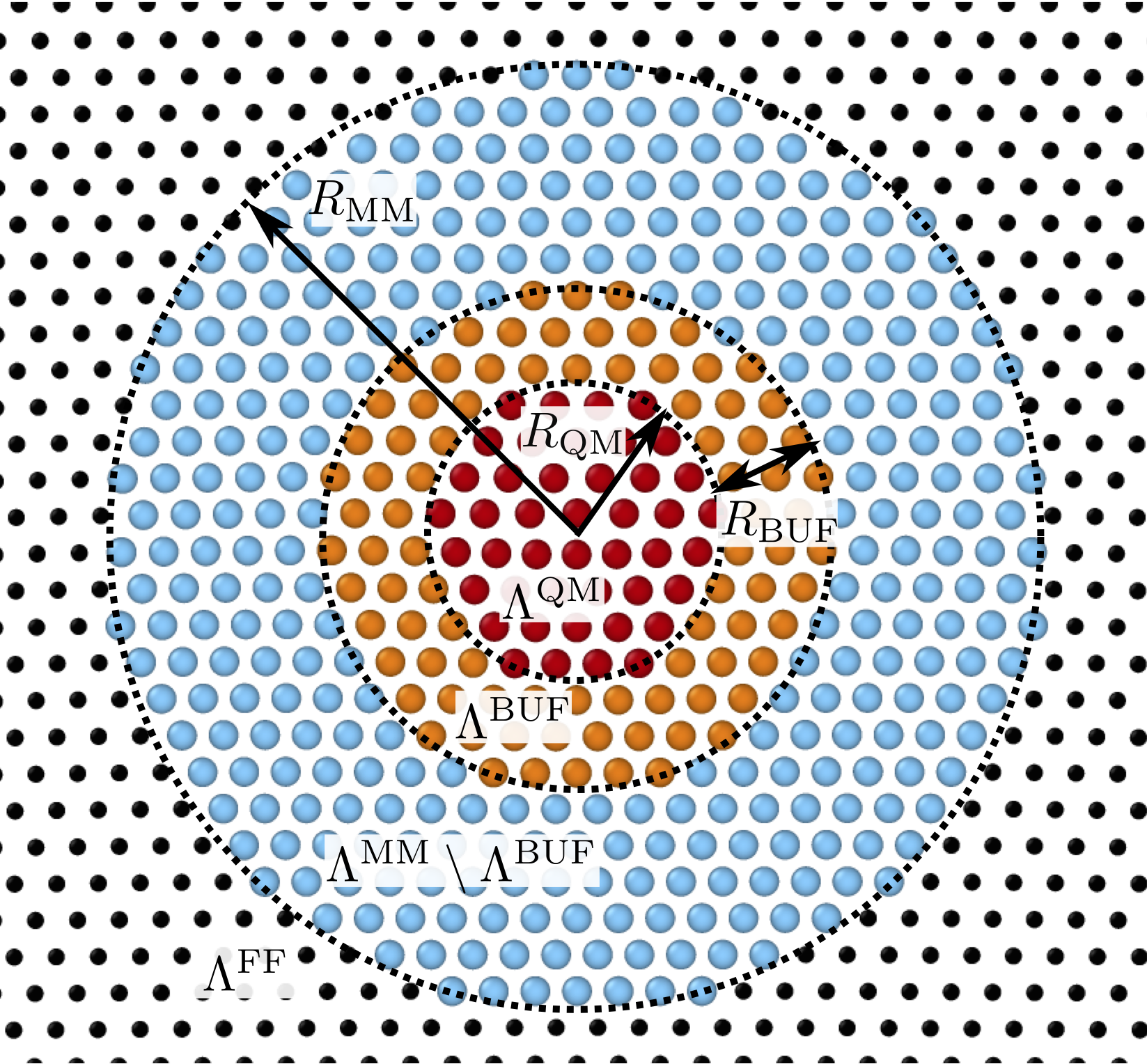}
        \caption{Decomposition of a crystal lattice into QM, MM, buffer, and far-field regions (the picture is from \cite[Figure 2]{chen17}).
		}
        \label{fig:decomposition}
\end{figure}

In QM/MM simulations, the reference configuration $\L$ is partitioned into three disjoint sets $\L = \LQM\cup \LMM\cup \LFF$, where $\LQM$ denotes the QM region containing the defect core, $\LMM$ denotes the MM region, and $\LFF$ denotes the far-field region with atomic positions frozen in the simulations. In addition to reduce spurious interface effects, we define a buffer region $\Lbuf\subset\LMM$ surrounding $\LQM$ such that all atoms in $\Lbuf\cup\LQM$ are involved in the QM calculations (to evaluate the site potentials in $\LQM$) but do not directly contribute to the hybrid energies of forces. 
For the sake of simplicity, we use balls centred at the defect core to decompose the domain $\Lambda$ with $\RQM$, $\RMM$ and $\Rbuf$, respectively (cf.~Figure~\ref{fig:decomposition}), but all methods and results have immediate generalizations to more complex domain decompositions.

\subsection{Accuracy of the interatomic potentials}
\label{sec:measurements}
The QM/MM coupling schemes can give good equilibrium state approximations only if, in the MM region, the MM site potential provides an accurate approximation to the QM site potential.
The decay of the displacement of equilibrium state \eqref{eq:ubar-decay} implies that, far from the defect core, the Taylor expansion of the QM site potential on homogeneous lattice (see \eqref{eq:taylor_E} or \cite{chen17}) will provide an excellent approximation.
%
This observation motivates us to introduce matching conditions between the QM and MM models in terms of Taylor expansion coefficients for the QM and MM potentials on the homogeneous lattice.

Let $\Vhom$ and $\VMM$ be the QM and MM homogeneous site potentials respectively, and assume both satisfy \assERL. Then we define the $j$th order site potential error by 
\begin{eqnarray}
\label{training_error}
\efit_j := \big|\delta^{j} V^{\rm h}(\pmb 0) - \delta^{j} \VMM(\pmb 0) \big|_{w^{-1}_j}
:= \left(\sum_{\pmb\rho = (\rho_1, ..., \rho_j)\in (\LhomS)^j}\Big|V^{\rm h}_{,\pmb\rho}(\pmb 0) - \VMM_{,\pmb\rho}(\pmb 0)\Big|^2 
w_j^{-1}(\pmb\rho) \right)^{1/2}\hspace{-1cm}
\end{eqnarray}
for $j=1,\cdots,\mathfrak{n}$,
where $V^{\rm h}_{,\pmb\rho}(\pmb 0)$ and $\VMM_{,\pmb\rho}(\pmb 0)$ are defined by \eqref{eq:dV}, and for $\gamma>0$,
$w_j:(\LhomS)^j \rightarrow \R$ is a weight function given by
\begin{align}
\label{eq:weight-def}
    w_j(\pmb \rho) := \prod_{i = 1}^j \wf(|\rho_i|) \quad \text{with}\quad \wf(|\rho|) = e^{-2\gamma|\rho|}. 
\end{align}
We do not have to introduce $\efit_0$ since adding a constant to the site potentials will not affect the equilibrium state, nor the energy differences. 
The choice of $\gamma$ only has a mild effect on the constants in our analysis. In practice, in constructing the loss function (see \eqref{cost:energymix} and \eqref{cost:forcemix}) of training MM site potential, we determine $\gamma$ in a way that roughly matches the locality of the site potentials \assERL.

Note that there is always some interaction range cutoff for an MM model, such that the interactions between two atoms with distance larger than the cutoff radius $R_{\rm cut}$ is ignored within the MM calculations.  
The choice of $R_{\rm cut}$ is important to determine the accuracy and computational cost of the MM model.  
The error from the neglected interactions is automatically included in $\efit_j$ since $\VMM_{,\pmb{\rho}}(\bm{0}) = 0$ when $|\rho_i|>R_{\rm cut}$ in \eqref{training_error}.

Alternatively, we can consider measurements on the accuracy of the forces.
We introduce the following error indicators comparing the $j$th order derivatives of the MM force $\FMM$ and QM force $\F^{\rm h}$ on the homogeneous lattice
\begin{eqnarray}
\label{training_error_F}
\ffit_j := \big| \delta^j\F^{\rm h}(\pmb 0) - \delta^j\FMM(\pmb 0)  \big|_{w^{-1}_j}
:= \left(\sum_{\pmb\ell = (\ell_1, ..., \ell_j)\in (\Lhom)^j}\Big|\F^{\rm h}_{,\pmb\ell}(\pmb{0}) - \FMM_{,\pmb\ell}(\pmb{0})\Big|^2
w_j^{-1}(\pmb\ell) \right)^{1/2} \hspace{-1cm}
\end{eqnarray}
for $j=1,\cdots,\mathfrak{n}-1$,
where $\mathcal{F}^{\rm h}_{,\pmb \ell}(\pmb{0})$ and $\mathcal{F}^{\rm MM}_{,\pmb \ell}(\pmb{0})$ are the $j$-th order derivatives of the QM and MM forces respectively, and $w_j$ is given by \eqref{eq:weight-def}.
Note that the forces vanish on the homogeneous lattice $\Lhom$ for both QM and MM models due to point symmetry, and therefore $\ffit_0$ is not required.
Similarly to the discussions for site potentials, the error from the truncation of interaction range at cutoff $R_{\rm cut}$ is again implicitly included.

We will show below that fitting an MM potential such that $\efit_j$ or $\ffit_j$ are small up to some order leads can lead to excellent QM/MM models. However, ensuring a good fit for high matching orders turns out to be very challenging. We will therefore propose also an alternative construction that sacrifices some accurary on $\efit_j$ or $\ffit_j$ but instead demands accurate fits of virials. The motivation for this is that this would incorporate continuum nonlinear elasticity which has been shown to provide an excellent model for atomistic interaction even fairly close to defect cores~\cite{weinan2007cauchy, co2013}. 


Thus, we briefly review the {\it Cauchy-Born rule}~\cite{weinan2007cauchy, co2013}, which relates the movement of atoms in a crystal to the overall deformation of the bulk solid.
For $\mathsf{F} \in \R^{d\times d}$, the Cauchy-Born rule makes an approximation such that in a crystalline solid subject to a small strain, the positions of the atoms within the crystal lattice follow the overall strain of the medium.
More precisely, the QM and MM site potentials are approximated by the Cauchy-Born elastic energy density functional $W^{\rm h}_{\rm cb}, ~ W^{\rm MM}_{\rm cb} : \R^{d\times d} \rightarrow \R$ respectively, with
\begin{eqnarray}
\label{eq:cbW}
W^{\rm h}_{\rm cb}\big(\mathsf{F}\big) := V^{\rm h}\big((\mathsf{F}-\mathsf{I})\cdot\LhomS\big) \qquad \text{and} \qquad W^{\rm MM}_{\rm cb}\big(\mathsf{F}\big) := V^{\rm MM}\big((\mathsf{F}-\mathsf{I})\cdot\LhomS\big).
\end{eqnarray}
The derivative (virial stress) and even higher order derivatives with respect to the deformation $\mathsf{F}$ can be obtained by direct calculations, 
\begin{eqnarray}
\label{dFWcb}
\partial^j_{\mathsf{F}} W_{\rm cb}^{\rm h}(\mathsf{F}_0) :=
\partial^j_{\mathsf{F}} W_{\rm cb}^{\rm h}(\mathsf{F})\big|_{\mathsf{F} = \mathsf{F}_0} = - \sum_{\pmb\rho=(\rho_1, \cdots, \rho_j) \in (\LhomS)^j} \Vhom_{,\pmb\rho} \big( (\mathsf{F}_0-\mathsf{I})\cdot\LhomS \big) \underbrace{\otimes{\rho}_1\otimes\cdots\otimes{\rho}_j}_{=:~\otimes\pmb\rho} ,
\end{eqnarray}
where $\mathsf{I}\in \R^{d\times d}$ is the identity matrix and $\otimes$ denotes the standard kronecker product.
We then introduce the following accuracy measure of the virial stress:
\begin{align}
\label{training_V}
\vfit_j := \big|\partial^{j+1}_{\mathsf{F}} \Wcb^{\rm h}(\mathsf{I}) - \partial^{j+1}_{\mathsf{F}} \Wcb^{\rm MM}(\mathsf{I}) \big|_{{w^{-1}_{j+1}}} :=\left(\sum_{\pmb\rho \in (\LhomS)^{j+1}}\Big|\big(V^{\rm h}_{,\pmb\rho}(\pmb 0) - \VMM_{,\pmb\rho}(\pmb 0)\big) \otimes {\pmb\rho}\Big|^2 w_{j+1}^{-1}(\pmb\rho) \right)^{1/2}
\end{align}
for $j=1,\cdots,\mathfrak{n}-1$, where the weights $w_j$ is defined by \eqref{eq:weight-def}.

\subsection{Energy-mixing}
\label{sec:energyMix}
The MM site potential $\VMM: (\R^d)^{\LhomS} \rightarrow \R$ is defined on the homogeneous reference configuration $\LhomS$, such that it can (i) be constructed uniformly and (ii) be compared with the QM site potential $\Vhom$ as shown in the previous subsection.
We will therefore need an interpolation operator to map the displacements on the defective reference configuration $\L$ to that on the corresponding homogeneous lattice $\Lhom$.
It is shown in \cite[Lemma D.1]{chen19} that such an interpolation operator $I^{\rm h}: \UsH(\L) \rightarrow \UsH(\Lhom)$ exists such that, for some constant $C_{*}>0$,
\begin{eqnarray}
\label{eq:eqI}
\|DI^{\rm h}u\|_{\ell^2_{\mathcal{N}}(\Lhom)}\leq C_{*}\|Du\|_{\ell^2_{\mathcal{N}}(\L)} \qquad \forall ~u \in \UsH(\L).
\end{eqnarray}
We emphasize that this interpolation is only introduced for the purpose of analysis, but does not affect the construction of the MM site potential or the coupling scheme, since the MM potential always has a cutoff and does not directly interact with the defect core region.

We now define the QM/MM hybrid energy difference functional by 
\begin{align}
\label{eq:hybrid_mid_energy}
\EH(u)  :=& \sum_{\ell\in \LQM} \Big( \VQM_{\ell}\big(Du_0(\ell) + D u(\ell)\big)  - \VQM_{\ell}\big(Du_0(\ell)\big) \Big) \nonumber \\
&+ \sum_{\ell\in \LMM\cup\LFF}  \Big( \VMM\big(Du_0(\ell) + D I^{\rm h}u(\ell)\big) - \VMM\big(Du_0(\ell)\big) \Big),
\end{align}
with $\VQM_{\ell}({\bf g}) := V_{\ell}^{\LQM \cup \Lbuf}({\bf g})$.

This energy difference functional was proposed in our previous paper \cite{chen17} (the second part of this series), with {\em strict} matching conditions, for $1 \leq j \leq \n-1$,
\begin{align}
\efit_j = 0, \quad \text{for energy-mixing} \qquad \text{or} \qquad \ffit_j = 0, \quad \text{for force-mixing},
\end{align}
which led to a convergent QM/MM scheme. However, in the present work we introduce {\em soft} matching conditions imposed via MLIPs and this will lead to a slight mismatch at the QM/MM interface (with error $\efit_1$ and $\efit_2$ respectively) that is difficult to converge to zero. 
This mismatch gives rise to the so-called ``ghost-force'' inconsistency; see Remark \ref{re:gfc} for the detailed mathematical reasoning. 
To reduce the effect of these ``ghost forces'' we employ ideas from the atomistic-to-continuum coupling literature~\cite{Chen12ghost, mlco2013, colz2012, co2011}.
%
%
Specifically, we use a ghost-force correction (GFC) \cite{colz2016}: to enforce the requirement that the homogeneous reference crystal is in equilibrium we make a dead load correction to the hybrid energy difference functional, 
\begin{eqnarray}
\label{eq:hybrid_energy_BGFC}
\quad \Egfc(u) 
:= \EH(u) - \<\delta\EH(\pmb{0}), \beta u\>, \qquad
\end{eqnarray}
where $\beta: \R^d \rightarrow \R$ is a characteristic function satisfying
\begin{eqnarray}
\label{beta}
	\beta(x) =
	\left\{ \begin{array}{ll}
		0, & |x|<\RQM/2
		\\[1ex]
		1, & |x|\geq\RQM/2  
	\end{array} \right. .
\end{eqnarray}

Note that $\beta$ can also be chosen as a smooth cutoff function that transiting from 0 to 1 as going away from the defect core (see \cite{mlcovkb2013}). 
If $\beta \in C^{2,1}(\R)$ satisfies the assumptions in \cite{colz2016}, then the error estimates in Theorem \ref{theo:energyM} will be improved due to the additional term $\|\nabla^2 \beta\|_{L^{\infty}}$ (see \cite[Section 5]{colz2016} for more details). However, in the present work, our aim is to construct a consistent energy-mixing scheme without ghost-forces, so the step function is sufficient for our analysis and will be used throughout this paper. The sharper estimates are thus not considered since the inefficiency of the energy-mixing schemes with GFC will also be shown in practice. 

The approximate equilibrium state of the defective system is then obtained by minimizing the energy difference functional with ghost force correction
\begin{eqnarray}
\label{eq:variational-problem-H}
\bar u^{\rm H} \in \arg\min \big\{ \mathcal{E}^{\rm GFC}(u), u \in \Admu^{\rm H}(\L) \big\} , 
\end{eqnarray}
where the admissible set is given by
\begin{eqnarray}
\label{e-mix-space}
\Admu^{\rm H}(\L):= \big\{ u\in\Us^{1,2}(\L) ~:~ y_0+u\in\Adm_0(\L) ~~{\rm and}~~u=0~{\rm in}~\Lambda^{\rm FF} \big\}.
\end{eqnarray}

The following result characterises previsely how the matching conditions $\efit_j$ affect the simulation error.
%
The proofs are given in \ref{sec:appendixProofEmix}.

\begin{theorem}
\label{theo:energyM}
	Suppose that \assERL~and either \asP~or \asD~are satisfied and that $\bar{u}$ is a strongly stable solution of \eqref{eq:variational-problem} satisfying \asS. Then, for any $K_{\rm E}\in\N,~2 \leq K_{\rm E} \leq \n-2$, $\RQM$ sufficiently large and $\efit_2$ sufficiently small, there exists an equilibrium $\uH$ of \eqref{eq:variational-problem-H} such that
	\begin{eqnarray}
	\label{res-bound}
	\|D\bar{u}-D\uH\|_{\ell^2_{\mathcal{N}}} 
	\leq C\left(\sum^{K_{\rm E}}_{j=2} \efit_{j}  R^{-\alpha_j}_{\rm QM} +   R^{-\alpha_{K_{\rm E}+1}}_{\rm QM} + R^{-d/2}_{\rm MM}\log^t\RMM + e^{-\kappa \Rbuf}
	\right) ,
	\end{eqnarray}
	where $C$ is a constant independent of $\RQM, \RMM, \efit$, and $\alpha_j$ and $t$ are given by
	\begin{align*}
	\alpha_j =
	\left\{\begin{aligned}
		&(2j-3)d/2 \quad 
		&\text{if}~\asP \\ 
		&j-2  
		&\text{if}~\asD
	\end{aligned}
	\right.
	\qquad{\rm and}\qquad
	t =
	\left\{\begin{aligned}
		&0 \quad 
		&\textrm{if}~\asP \\ 
		&1 
		&\textrm{if}~\asD
	\end{aligned}
	\right..
	\end{align*}
\end{theorem}

Our second result on energy mixing schemes demonstrates how adding soft matching conditions for the virials improves the rate in terms of the QM region size. As we will see below, this comes at essentially no additional computational cost.

\begin{theorem}
\label{thm:egyM_V}
Under the conditions of Theorem \ref{theo:energyM}, for any $ 2 \leq K_{\rm E} \leq \n-2$, we have
	\begin{align}
	\label{res-bound-E-stress}
	\|D\bar{u}-D\uH\|_{\ell^2_{\mathcal{N}}} 
	\leq C\left( \sum^{K_{\rm E}}_{j=2} \efit_j  R^{-\alpha_j}_{\rm QM} + \vfit_{K_{\rm E}}  R^{-\beta_{K_{\rm E}}}_{\rm QM} +  R^{-\beta_{K_{\rm E}}-1}_{\rm QM} + R^{-d/2}_{\rm MM}\log^t\RMM + e^{-\kappa \Rbuf}\right)
	\end{align}
	with
	\begin{align*}
	\alpha_j =
	\left\{\begin{aligned}
		&(2j-3)d/2 ~~ 
		&\textrm{if}~\asP \\ 
		&j-2
		&\textrm{if}~\asD
	\end{aligned} 
	\right.,
	\quad
	\beta_{K_{\rm E}} =
	\left\{\begin{aligned}
		&\alpha_{K_{\rm E}}+d ~~ 
		& \textrm{if}~\asP \\ 
		&\alpha_{K_{\rm E}}+1
		&\textrm{if}~\asD
	\end{aligned}
	\right.
	\quad{\rm and}\quad
	t =
	\left\{\begin{aligned}
		&0 ~~
		&\textrm{if}~\asP \\ 
		&1
		&\textrm{if}~\asD
	\end{aligned}
	\right..
\end{align*} 
\end{theorem}

The error estimates in the foregoing theorems identify how the error in the equilibrium geometry depends on they key approximation parameters, model accuracy, QM region size and MM region size. In \cite[Theorem 4.1]{chen17}, it was assumed that $\efit_j=0$ for $j=1, \cdots, K_{\rm E}$ which then yields immediate rates of convergence in terms of QM and MM region size. For both theorems, one can then readily balance the $\RQM$ and $\RMM$ parameters to obtain a single rate. These optimal balanced are summarized in Table~\ref{table-e-mix}. Even with $\efit_j=0$, Theorem~\ref{thm:egyM_V} is a new result. Here, a key observation that will lead to particularly competitive QM/MM schemes is that the rate of convergence improves by a full power when $\vfit_{K_E} = 0$. 

In the present work, our focus is the dependence of the estimates on the $\efit_j, \ffit_j$ and $\vfit_j$ parameters. By fitting MM potentials such that they are non-zero but ``sufficiently small'' we will obtain models where the rates from Table~\ref{table-e-mix} appear in the pre-asymptotic regime. This will be clearly observed in the numerical tests of Section~\ref{sec:numerics}.

\subsection{Force-mixing}
\label{sec:forceMix}
%
To avoid the complications surrounding the ghost-forces correction (GFC), force-mixing schemes design the coupling in terms of QM and MM forces instead of energies. 
The main advantage of force-mixing is that there are no spurious interface forces as in energy-mixing. Even though GFC can be applied in energy-mixing, the force-mixing schemes are still considered more convenient and efficient than energy-mixing in practical simulations \cite{csanyi04, kermode08}.

The approximate equilibrium state is then obtained by solving the following hybrid force balance equations: {\it Find $\uH \in \Admu^{\rm H}(\L)$ such that}
\begin{eqnarray}
\label{problem-f-mix}
	\F_{\ell}^{\rm H}\big(\bar{u}^{\rm H}\big) = 0
	\qquad \forall~\ell\in\LQM\cup\LMM
\end{eqnarray}
where $\F_{\ell}^{\rm H}$ are the hybrid forces
\begin{eqnarray}
\label{F-H}
	\F_{\ell}^{\rm H}(u) =
	\left\{ \begin{array}{ll}
		\F^{\rm QM}_{\ell}(u), \quad & \ell\in\LQM
		\\[1ex]
		\FMM_\ell(u),
		 \quad & \ell\in\LMM  
	\end{array} \right. ,
\end{eqnarray}
with $\FMM_\ell(u):=\FMM\big( u_0(\cdot-\ell) + I^{\rm h}u(\cdot-\ell)\big)$ and
$I^{\rm h}: \UsH(\L) \rightarrow \UsH(\Lhom)$ the interpolation given in Section \ref{sec:energyMix} and $\F^{\rm QM}_{\ell}(u) := \F^{\LQM \cup \Lbuf}_{\ell}(u)$. 
The error estimates for the QM/MM force-mixing scheme \eqref{problem-f-mix} are stated in the following theorems, whose proofs are again postponed to \ref{sec:appendixProofFmix}.

\begin{theorem}
	\label{theo:forceM}
	Suppose that \assERL~and either assumption \asP~or \asD~are satisfied and that $\bar{u}$ is a strongly stable solution of \eqref{eq:problem-force} satisfying \asS. Then, for any $K_{\rm F}\in\N,~1 \leq K_{\rm F} \leq \n-3$, $\RQM$ sufficiently large satisfying $\displaystyle \log\frac{\RMM}{\RQM} \leq C^{\rm MM}_{\rm QM}$, and $\ffit_1$ sufficiently small, there exists an equilibrium $\uH$ of \eqref{problem-f-mix} such that 
	\begin{align}\label{res-bound-F}
	\|D\bar{u}-D\uH\|_{\ell^2_{\mathcal{N}}} \leq C\left( \Big( \sum^{K_{\rm F}}_{j=1} \ffit_j R^{-\alpha_j}_{\rm QM} + R^{-\alpha_{K_{\rm F}+1}}_{\rm QM} \Big)\log\RMM + R^{-d/2}_{\rm MM}\log^t\RMM+e^{-\kappa \Rbuf}\right),
	\end{align}
	where $C$ is a constant independent of $\RQM, \RMM, \efit$, and $\alpha_j$ and $t$ are given by
		\begin{align*}
	\alpha_j =
	\left\{\begin{aligned}
		&(2j-1)d/2 \quad 
		&\text{if}~\asP \\ 
		&j-1  
		&\text{if}~\asD
	\end{aligned}
	\right.
	\quad{\rm and}\quad
	t =
	\left\{\begin{aligned}
		&0 \quad 
		&\textrm{if}~\asP \\ 
		&1 
		&\textrm{if}~\asD
	\end{aligned}
	\right..
	\end{align*}
\end{theorem}
	
\begin{theorem}
\label{thm:forceM_V}
	Under the conditions of Theorem \ref{theo:forceM}, for any $1\leq K_{\rm F} \leq \n-3$, we have
	\begin{align}\label{res-bound-F-stress}
	\|D\bar{u}-D\uH\|_{\ell^2_{\mathcal{N}}} \leq C\Bigg( \Big( \sum^{K_{\rm F}}_{j=1} \ffit_j  R^{-\alpha_j}_{\rm QM} +& \vfit_{K_{\rm F}+1} R^{-\beta_{K_{\rm F}}}_{\rm QM} +  R^{-\beta_{K_{\rm F}}-1}_{\rm QM} \Big)\log\RMM
	\nonumber \\
	&+  R^{-d/2}_{\rm MM}\log^t\RMM + e^{-\kappa \Rbuf} \Bigg)
	\end{align}
	with
\begin{align*}
	\alpha_j =
	\left\{\begin{aligned}
		&(2j-1)d/2 \quad 
		&\textrm{if}~\asP \\ 
		&j-1
		&\textrm{if}~\asD
	\end{aligned} 
	\right.,
	\quad
	\beta_{K_{\rm F}} =
	\left\{\begin{aligned}
		&\alpha_{K_{\rm F}}+d \quad 
		& \textrm{if}~\asP \\ 
		&\alpha_{K_{\rm F}}+1
		&\textrm{if}~\asD
	\end{aligned}
	\right.
	\quad{\rm and}\quad
	t =
	\left\{\begin{aligned}
		&0 \quad
		&\textrm{if}~\asP \\ 
		&1
		&\textrm{if}~\asD
	\end{aligned}
	\right..
\end{align*} 
\end{theorem}

\begin{remark}\label{log}
	According to the assumption $\log \frac{\RMM}{\RQM} \leq C^{\rm MM}_{\rm QM}$, $\RMM$ is in fact bounded in terms of a polynomial of $\RQM$, and we could simply replace $\log\RMM$ with $\log\RQM$ in \eqref{res-bound-F} and \eqref{res-bound-F-stress}. However, we keep $\log\RMM$ to highlight the dependence of the error estimate on the growth of $\RMM$ relative to $\RQM$, and thus explain the necessity of carefully controlling the $\RMM$ parameter; see also \cite[Remark 5.1]{chen17}. 
\end{remark}

The above estimates for the force-mixing scheme are consistent with our previous result in \cite[Theorem 5.1]{chen17}. The optimal balance of $\RQM$ and $\RMM$ and resulting estimates for the force-mixing schemes are again summarized in Table \ref{table-e-mix}, under the assumption that all matching parameters satisfy $\ffit_j = 0$. 
%
As in the energy-mixing schemes we observe that matching the $(K_{\rm F}+1)$-th derivative of the virial stress significantly improves the convergence rates of the coupling scheme.
Moreover, we observe that the force-mixing schemes can be much cheaper than the energy-mixing schemes to obtain the same convergence rate, since lower order derivatives need to be matched with the QM model for force-mixing. For example, to achieve the $(2m+1)$-th order convergence for \asP~in two dimension, $O(n^m)$ observations are needed for force-mixing while $O(n^{m+1})$ for energy-mixing, where $n$ is the number of the neighbours of each site. The MLIPs for force-mixing are then fitted by using less approximation parameters to achieve the same accuracy. Hence, the force-mixing schemes are of particular practical interest. See Table \ref{table-params} and the discussions in Section \ref{sec:simulations} for more quantitative details.



\begin{table}
	\begin{center}
    {\bf Point defects in two dimension}\\
	\vskip0.2cm
	\begin{tabular}{|c|ccc|ccc|}
	\hline
		{Case {\bf P}, $d=2$}& \multicolumn{3}{c|}{no conditions on $\vfit_K$} &
		\multicolumn{3}{c|}{$\displaystyle \vfit_{K} \leq \epsilon$ ($K=K_{\rm E}$ or $K_{\rm F}+1$)} \\[1mm]
		\hline
		$K_{\rm E}$  &  2  & 3 & 4 & 2 & 3 & 4  \\
		$K_{\rm F}$  &  1  & 2 & 3 & 1 & 2 & 3  \\
		\hline
		$\RMM$    &  $\RQM^3$ & $\RQM^5$ & $\RQM^7$
		& $\RQM^4$ & $\RQM^6$ & $\RQM^8$
	    \\[1mm]
	    \hline
		error     &  $\RQM^{-3}$ & $\RQM^{-5}$ & $\RQM^{-7}$
		& $\RQM^{-4}$ & $\RQM^{-6}$ & $ \RQM^{-8}$
		 \\[1mm]
		\hline
	\end{tabular}
	\vskip 1.3cm
	{\bf Point defects in three dimension}\\
	\vskip0.2cm
	\begin{tabular}{|c|ccc|ccc|}
	\hline
		Case {\bf P}, $d=3$ & \multicolumn{3}{c|} {no conditions on $\vfit_K$}&
		\multicolumn{3}{c|}{$\displaystyle \vfit_{K} \leq \epsilon$ ($K=K_{\rm E}$ or $K_{\rm F}+1$)} \\[1mm]
		\hline
		$K_{\rm E}$ & 2 & 3 & 4 & 2 & 3 & 4 \\
		$K_{\rm F}$ & 1 & 2 & 3 & 1 & 2 & 3  \\
		\hline
		$\RMM$ & $\RQM^3$ & $\RQM^5$ & $\RQM^7$
		 & $\RQM^{11/3}$ & $\RQM^{17/3}$ & $\RQM^{23/3}$
		 \\[1mm]
		 \hline
		error & $\RQM^{-4.5}$ & $\RQM^{-7.5}$ & $\RQM^{-10.5}$
		& $\RQM^{-5.5}$ & $\RQM^{-8.5}$ & $ \RQM^{-11.5}$ \\[1mm]
		\hline
	\end{tabular}
	\vskip 1.3cm
	{\bf Straight Dislocations}\\
	\vskip0.2cm
	\begin{tabular}{|c|ccc|ccc|}
	\hline
		Case {\bf D} & \multicolumn{3}{c|} {no conditions on $\vfit_K$}&
		\multicolumn{3}{c|}{$\displaystyle \vfit_{K} \leq \epsilon$ ($K=K_{\rm E}$ or $K_{\rm F}+1$)} \\[1mm]
		\hline
		$K_{\rm E}$  & 2 & 3 & 4 & 2 & 3 & 4 \\
		$K_{\rm F}$  & 1 & 2 & 3 & 1 & 2 & 3  \\
		\hline
		$\RMM$ &  $\RQM$ & $\RQM^2$ & $\RQM^3$
		 & $\RQM^2$ & $\RQM^3$ & $\RQM^4$
		 \\[1mm]
		 \hline
		error    &  $\RQM^{-1}$ & $\RQM^{-2}$ & $\RQM^{-3}$
		& $\RQM^{-2}$ & $\RQM^{-3}$ & $ \RQM^{-4}$ \\[1mm]
		\hline
	\end{tabular}
	\vskip 1.2cm
	\caption{The optimal choice of $\RMM$ and the error decay with respect to $\RQM$ for both point defects and dislocations, by using energy-mixing (varying $K_{\rm E}$) and force-mixing (varying $K_{\rm F}$) coupling schemes.
	}
	\label{table-e-mix}
	\end{center}
\end{table}


\section{Numerical experiments}
\label{sec:numerics}
\setcounter{equation}{0}
In this section, we present a concrete implementation of QM/MM schemes inspired by our theory, and use these to confirm our analytical error estimates on several model systems.

\subsection{Constructions of the MLIPs}
\label{sec:constructions}
\setcounter{equation}{0}
%
%
We will construct MLIPs to be used as the MM model in QM/MM coupling schemes. Following our theoretical results we will require that the MLIPs match a QM reference model in the sense that $\efit_j,~\ffit_j$ and $\vfit_j$ are small.
%
%

{\it Step 1 - Training Set: } Towards that end we first generate a training set: For energy-mixing schemes \eqref{eq:variational-problem-H} the training set contains the observations $\delta^j V^{\rm h}(\bm 0)$ for $j=1, 2, ..., K_{\rm E}$, where $K_{\rm E} \in \N$ determines the order of accuracy of the MM model (cf. Theorem~\ref{theo:energyM}). 
According to Theorem~\ref{thm:egyM_V} we can further improve the accuracy of the model by adding the virial stress observations $\partial^{K_{\rm E}+1}_{\mathsf{F}} W_{\rm cb}^{\rm h}(\mathsf{I})$ to the training set. In practise we will restrict all our tests to $K_E \leq 3$ and to $K_E \leq 2$ with virials.  

For force-mixing schemes \eqref{problem-f-mix}, the training set contains the derivatives of the QM force $\delta^j \mathcal{F}^{\rm h}({\bm 0})$ for $j=1, 2, ..., K_{\rm F}$ at the homogeneous lattice $\Lhom$, with a given $K_{\rm F} \in \N$, where $K_{\rm F}$ now determines the MM model accuracy. Again, we can add the $(K_{\rm F}+1)$-th order derivatives of the virial, $\partial^{K_{\rm F}+2}_{\mathsf{F}} W_{\rm cb}^{\rm h}(\mathsf{I})$ to the training set to improve the convergence rate; cf. Theorem \ref{thm:forceM_V}. For force-mixing we restrict $K_{\rm F} \leq 2$. 

When creating the training sets from the QM reference models, the zeroth and first order derivatives (of either site energies or forces) are calculated from first order perturbation theory of eigenvalue problems \cite{kato}, i.e., close to machine precision. Higher order derivatives of site energies, forces or virials are obtained from centered finite difference approximations, which provides at least half of machine precision and is more than enough for our purposes. Note also that these finite difference approximations are only applied during the training but not the evaluation of the MLIP.

{\it Step 2 - Parameterisation: } Next, we need to choose a parameterisation of the MM model to match the observations contained in the training set. Since we require a fairly stringent accuracy on a wide training set we will use an MLIP ansatz. Although a wide variety of choices is available nowadays, we have opted for the atomic cluster expansion (ACE) \cite{bachmayr19,oord19,2020-pace1} which is a {\em linear} model but still outperforms most alternative approaches~\cite{2020-pace1}. We do not believe that the specific choice of MLIP is essential, though, and in fact expect that most MLIPs would lead to similar results with suitable minor  modifications of our procedures. For this reason, we do not present the details in the main text, but provide a review of the ACE model in \ref{sec:ACE}. For the purpose of setting up the MM potential and QM/MM method we only require that the ACE model provides a parameterised site potential (see also \eqref{ships:energy}), for $\pmb{g}\in\big(\R^d\big)^{\LhomS}$, 
\[
V^{\rm ACE}(\pmb{g}; \{c_B\}_{B\in\pmb{B}}) = \sum_{B\in\pmb{B}} c_B B(\pmb{g}) 
\]
where $B$ are the basis functions and $c_B$ the parameters. Given this parameterised site potential, we can evaluate site energies, total energies, forces, and virials.

{\it Step 3 - Parameter estimation: } Given the training set and parameterisation, we determine the parameters by minimising a suitable loss functional. For energy mixing schemes  we use the loss 
\begin{align}
\label{cost:energymix}
\mathcal{L}_{\rm E}\big(\{c_B\}\big) := \sum_{j=0}^{K_{\rm E}} W^{\rm E}_j\big(\varepsilon_j^{\rm E}\big)^2
= \sum_{j=0}^{K_{\rm E}}\left(W^{\rm E}_j\sum_{\pmb\rho = (\rho_1, ..., \rho_j)\in (\LhomS)^j}\Big|V^{\rm h}_{,\pmb{\rho}}(\pmb 0) - V^{\rm ACE}_{,\pmb\rho}(\pmb 0; \{c_B\})\Big|^2 {w}_j^{-1}(\pmb\rho) \right),
\end{align}
where $w_j(\pmb\rho)$ is the weight function defined by \eqref{eq:weight-def} and $W_j^{\rm E}$ are additional weights that may depend on the configurations and the observations (i.e. $V^{\rm h}_{,\pmb{\rho}}(\pmb 0)$) and are used to up- or down-weight specific observations. We determine the additional weights $W^{\rm E}_j$ and below also $W^{\rm F}_j$ and $W^{\rm V}_j$ purely empirically and provide the details in Table \ref{table-weights} in \ref{sec:ns}.

%
Similarly, for constructing a force-mixing scheme we minimise the loss 
\begin{align}
\label{cost:forcemix}
\mathcal{L}_{\rm F}\big(\{c_B\}\big) := \sum_{j=0}^{K_{\rm F}} W^{\rm F}_j \big(\varepsilon_j^{\rm F}\big)^2
= \sum_{j=0}^{K_{\rm F}} \left(W^{\rm F}_j\sum_{\pmb\ell = (\ell_1, ..., \ell_j)\in (\Lhom)^j}\Big|\F^{\rm h}_{,\pmb{\ell}}(\pmb 0) - \F^{\rm ACE}_{,\pmb{\ell}}(\pmb 0; \{c_B\})\Big|^2 {w}_j^{-1}(\pmb\ell) \right)
\end{align}
If derivatives of the virial stress are included in the training set, then the loss functions are further improved by adding the term 
\begin{align}
\label{cost:forcemix-V}
W^{\rm V}_K \big(\varepsilon_K^{\rm V} \big)^2.
\end{align}
%

In all cases, the loss function is quadratic in the parameters $(c_B)_{B\in \pmb{B}}$ and can therefore be minimised using a QR factorisation. In our implementation we use a rank-revealing QR (rr-QR) factorisation \cite{chan1987rank} which provides a mechanism analogous to Tychonov-regularisation~\cite{oord19}. The regularisation parameter for rr-QR factorisation is set to be $10^{-5}$ throughout our numerical experiments.

\subsection{Numerical results}
\label{sec:simulations}
We now present numerical tests for two prototypical defective crystals: a point defect and a straight edge dislocation. While the edge dislocation can naturally be simulated in a quasi-2D setup, we will also treat the point defect in an analogous quasi-2D manner, in order to control the computational cost. Thus all geometries we consider are taken infinite in the $(x, y)$-plane but with periodic boundary conditions in the $z$-direction.

We will perform these tests for two QM reference models: (1) First, in order to test our ideas in the simplest possible setting we will use an embedded atom model (EAM) \cite{Daw1984a} as the reference model, instead of an actual QM model. This allows us to explore our matching conditions to relatively high order and thus compare different schemes in a wider parameter range. (2) After identifying the most promising QM/MM candidate schemes, we will then test those on the more realistic NRL tight binding model \cite{cohen94,mehl96,papaconstantopoulos15} (see also \ref{sec:NRL} for a brief introduction), as the QM reference model.

All simulations are implemented in open-source {\tt Julia} packages {\tt SKTB.jl} \cite{gitSKTB} (for the NRL tight binding model), {\tt ACE.jl} \cite{gitACE} (for the ACE model) and {\tt QMMM2.jl} \cite{gitQMMM2} (for the QM/MM coupling scheme). All tests we report are performed on an {\tt Intel(R) Core(TM) i7-7820HQ CPU @2.90GHz}, with {\tt macOS (x86-64-apple-darwin19.6.0)} operating system.

\subsubsection{Ghost-forces and patch test}
Before we embark on convergence tests, we numerically illustrate the effect of ghost forces in the first energy mixing scheme \eqref{eq:hybrid_mid_energy},
and thus the need for the ghost force correction scheme \eqref{eq:hybrid_energy_BGFC}.

Following the instructions in Section \ref{sec:constructions}, we construct the MLIP by minimising the loss function \eqref{cost:energymix} with $K_{\rm E}=2$, where the weights are given in Table \ref{table-weights}. The total number of the observations and parameters $(c_B)_{B\in \pmb{B}}$ are, respectively, 6001 and 668. The MLIP is then constructed such that the relative root mean square error (RRMSE) of $\efit_1$ and $\efit_2$ are both less than $2\%$.



\begin{figure}[!htb]
	\centering 
	\subfigure[The ghost force in the system with a single vacancy.]{
		\label{fig:gfa}
		\includegraphics[height=6cm]{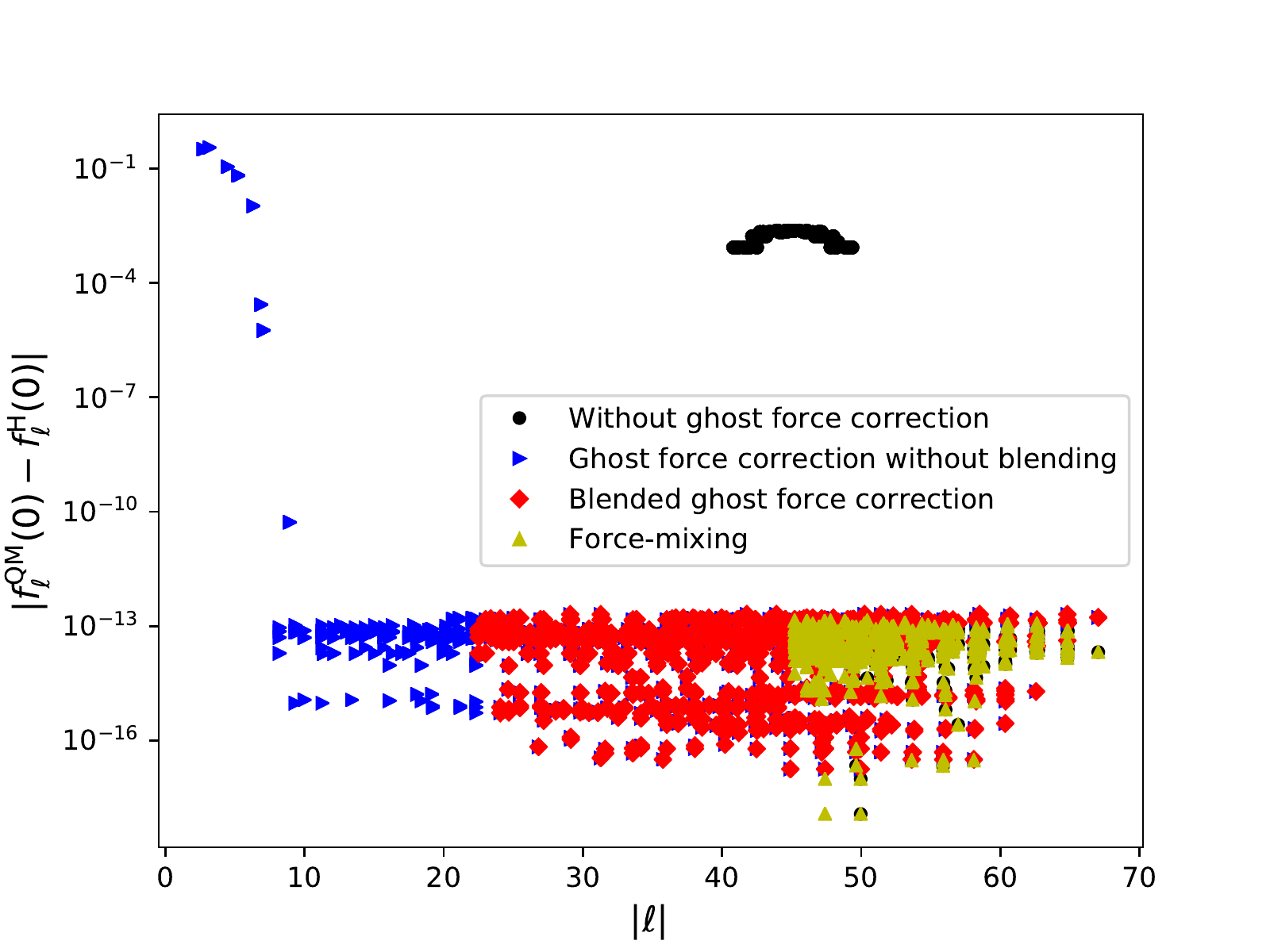}}
	\hspace{0.2cm} 
	\subfigure[Patch test in the homogeneous lattice.]{
		\label{fig:gfb} 
		\includegraphics[height=6cm]{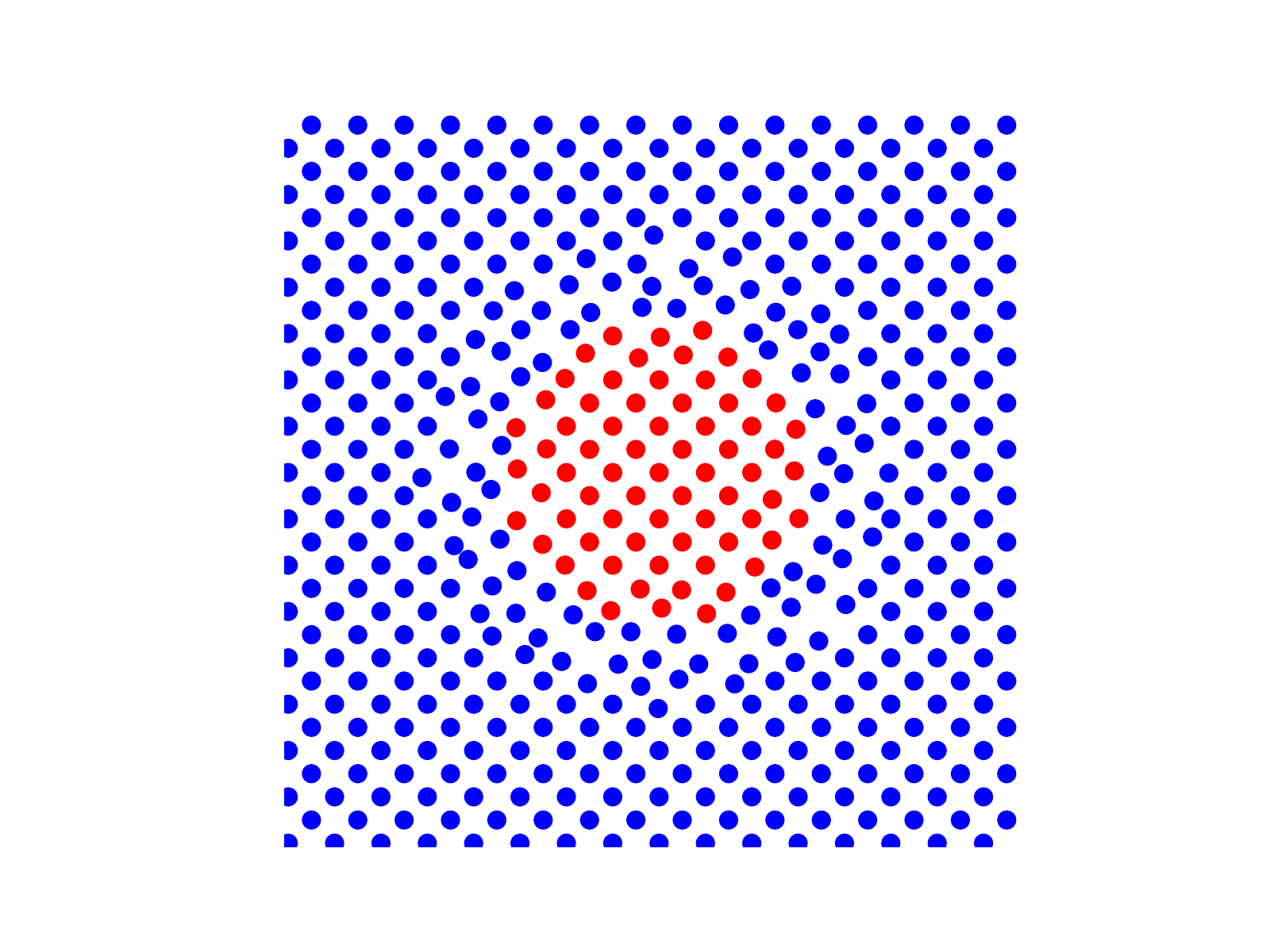}}
	\caption{Ghost force correction and patch test.}
	\label{fig:gf}
\end{figure}

In Figure \ref{fig:gf}(a) we plot the hybrid forces $f^{\rm H}_{\ell}(\p0)$ produced by the four QM/MM energy mixing scheme in an unrelaxed system containing a vacancy at the origin. We can observe the spurious forces (ghost forces) in the QM/MM interface region for the uncorrected energy mixing scheme. Of particular note is that, despite the fairly stringent accuracy we achieved for $\efit_1$, these forces are almost of the same magnitude as the forces in he defect core (before relaxation) and will therefore have a significant influence on the relaxation. For all other schemes, there are no such spurious interface forces. 

If we minimize the energy functional \eqref{eq:hybrid_mid_energy} (without ghost force correction) starting at a homogeneous lattice configuration we can observe the effect of the ghost forces on the equilibrium (patch test). A snapshot of the relaxation process is shown in Figure \ref{fig:gf}(b), where the red and blue points denote the atoms in QM and MM region respectively.  We observe that the homogeneous lattice configuration is unstable during the energy minimisation.


\subsubsection{Convergence results for the EAM W model}
Next, we explore a variety of QM/MM models where the reference model is given by an EAM potential \cite{Daw1984a} for W. Using an EAM potential as the reference model instead of an actual ab initio model means that we can more easily perform large-scale tests in a wider parameter regime in order to narrow down the best choices.

We produce both energy and force-mixing schemes, where the MLIPs are fitted by following the construction in Section \ref{sec:constructions}. The weights are given in Table \ref{table-weights}. The detailed fitting parameters, fitting accuracy in terms of the relative root mean square error (RRMSE) and fitting time are given in Table \ref{table-params}, where the total number of the observations and parameters $(c_B)_{B\in \pmb{B}}$ are denoted by $\# {\pmb O}$ and $\# \pmb{B}$, respectively. We also note that the $\efit_0$ and $\ffit_0$ are not required since they will not affect the equilibrium state.
 
\begin{table}
	\begin{center}
    {\bf Fitting accuracy on measurements and fitting time T for energy-mixing}\\
	\vskip0.2cm
	\begin{tabular}{|c|c|c|c|c|c|}
	\hline
		\diagbox{Schemes}{Accuracy (RRMSE)}{Measurements} & {$\efit_1$} & {$\efit_2$} & {$\efit_3$} & {$\vfit_2$} & {T~(s)}\\[0.5mm]
		\hline
		$K_{\rm E}=2$, $\# {\pmb O}=6001$, $\# {\pmb B}=668$   & 1.84\% & 0.42\% & - & - & 61.5  \\[0.5mm]
		\hline
		$K_{\rm E}=3$, $\# {\pmb O}=76201$, $\# {\pmb B}=5456$   & 5.46\% & 1.33\% & 1.09\% & - & 4323.9
	    \\[0.5mm]
	    \hline
		$K_{\rm E}=2 ~~\&$  virial, $\# {\pmb O}=6730$, $\# {\pmb B}=668$   & 3.90\% & 0.88\% & - & 0.45\% & 84.6
		\\[0.5mm]
		\hline
	\end{tabular}
	\\[6mm]
    {\bf Fitting accuracy on measurements and fitting time T for force-mixing}\\
	\vskip0.3cm
	\begin{tabular}{|c|c|c|c|c|c|}
	\hline
		\diagbox{Schemes}{Accuracy (RRMSE)}{Measurements} & {$\ffit_1$} & {$\ffit_2$} & {$\vfit_2$} & {$\vfit_3$} & {T (s)}\\[0.5mm]
		\hline
		$K_{\rm F}=1$, $\# {\pmb O}=600$, $\# {\pmb B}=110$ & 0.005\% & - & - & - & 1.8  \\[0.5mm]
		\hline
		$K_{\rm F}=2$, $\# {\pmb O}=6000$, $\# {\pmb B}=668$ & 0.29\% & 0.41\% & - & - & 40.6
	    \\[0.5mm]
	    \hline
		$K_{\rm F}=1~~\&$ virial, $\# {\pmb O}=1329$, $\# {\pmb B}=110$ & 0.58\% & - & 0.28\% & - & 5.8
		\\[0.5mm]
	   \hline
	     $K_{\rm F}=2~~\&$ virial, $\# {\pmb O}=12561$, $\# {\pmb B}=1119$ & 0.94\% & 1.07\% & - & 0.61\% & 348.3
		\\[0.5mm]
		\hline
	\end{tabular}
	\end{center}
	\vskip 0.5cm
	\caption{EAM W reference model: The fitting accuracy (RRMSE) on different measurements and fitting time for energy-mixing (varying $K_{\rm E}$) and force-mixing (varying $K_{\rm F}$) coupling schemes.
	}
	\label{table-params}
\end{table}

\begin{figure}[!htb]
    \centering
    \includegraphics[scale=0.48]{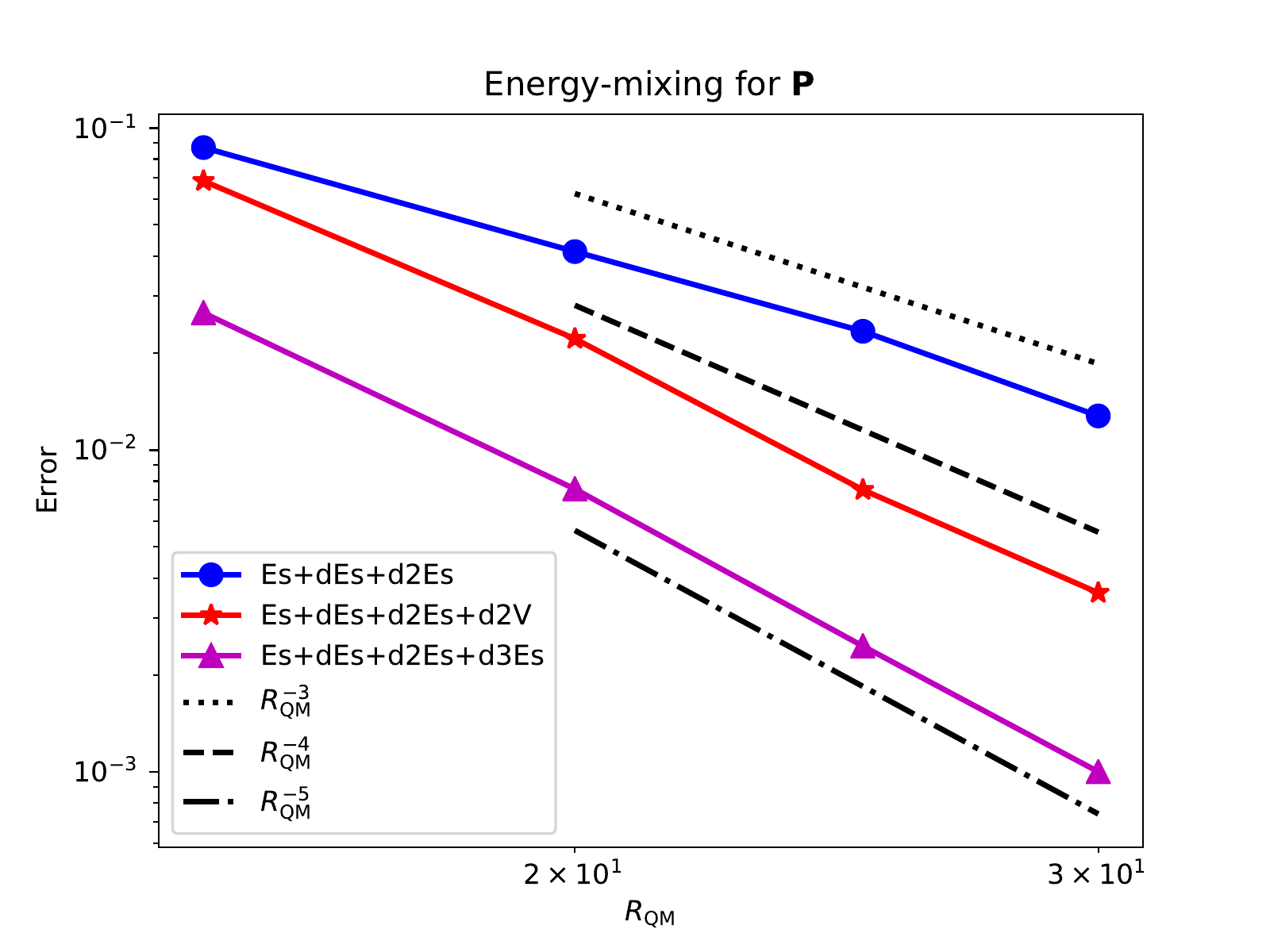}
    \includegraphics[scale=0.48]{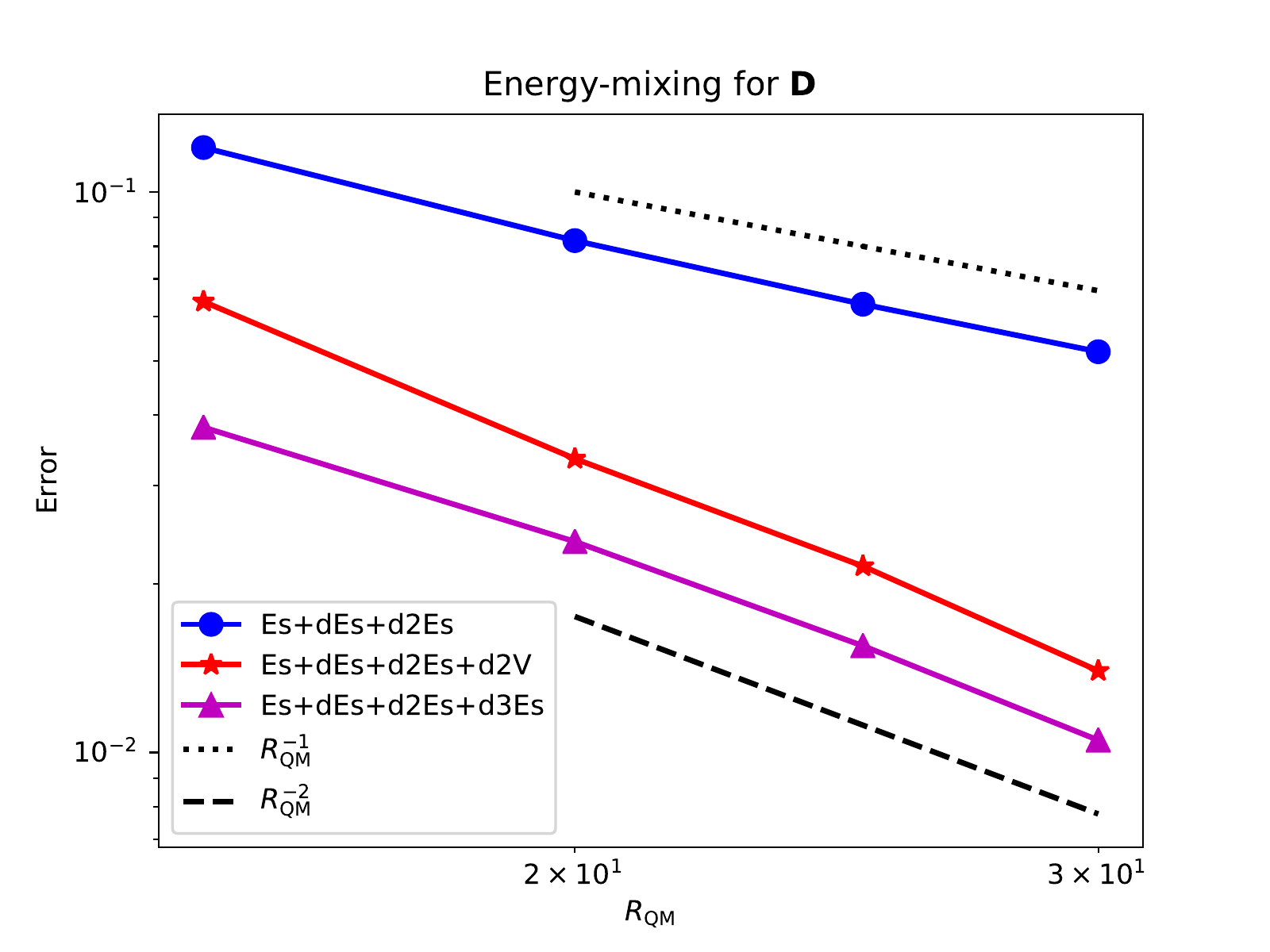}
    \caption{EAM W reference model: Convergence rates of energy-mixing scheme for point defects (left) and dislocations (right).
    }
    \label{fig:energy-mixing}
\end{figure}

\begin{figure}[!htb]
    \centering
    \includegraphics[scale=0.48]{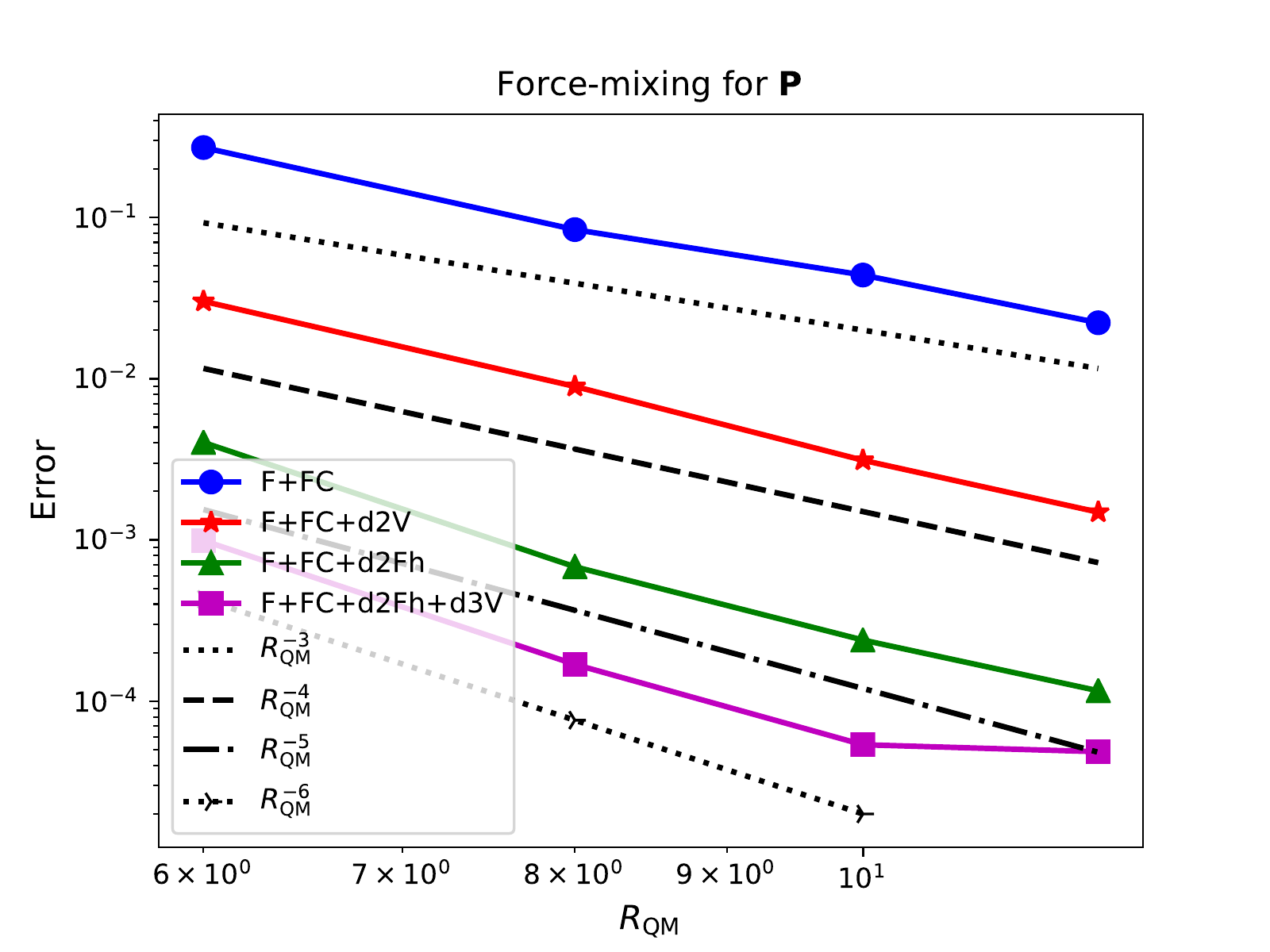}
    \includegraphics[scale=0.48]{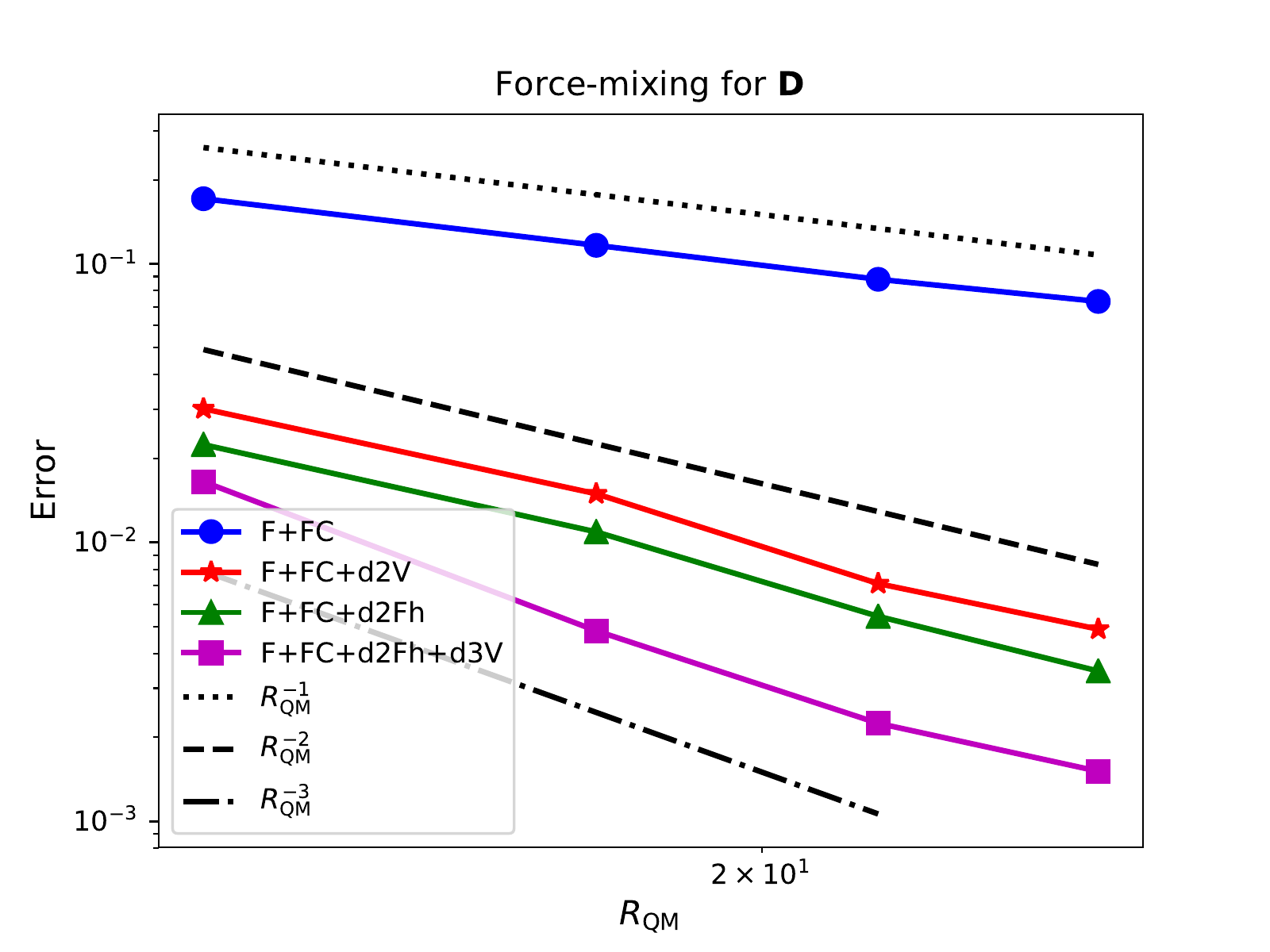}
    \caption{EAM W reference model: Convergence rates of the force-mixing scheme for point defects (left) and dislocations (right).}
    \label{fig:force-mixing}
\end{figure}

We then equilibrate our QM/MM models for (case \asP) a vacancy defect; and (case \asD) a (001)[100] straight edge-dislocation. 
The $\RMM$ and $\Rbuf$ parameters are chosen sufficiently large {($\RMM=120r_0$, $\Rbuf=6r_0$, where $r_0$ is the lattice constant of W)}, such that the error from the truncation of $\RQM$ will dominate. 
The decay of approximate error $\|D\bar{u} - D\bar{u}^{\rm H}\|_{\ell^2_{\mathcal{N}}}$ of the QM/MM coupling schemes are shown in Figure \ref{fig:energy-mixing} and \ref{fig:force-mixing} for, respectively, the energy-mixing and force-mixing schemes.
We observe that the convergence rates perfectly match our theoretical predictions from Theorems \ref{thm:forceM_V} and which summarized in Table \ref{table-e-mix}. 

Combing the results in Table \ref{table-params}, Figure \ref{fig:energy-mixing} and Figure \ref{fig:force-mixing}, we observe that the force-mixing schemes are of particular practical interest. MLIPs for force-mixing achieve much better fit accuracy at much lower number of parameters and lower evaluation cost (directly proportional to the number of parameters for ACE), since fewer observations are required for force-mixing to achieve the same accuracy. Our second important observation is that matching derivatives of the virial stress significantly improves the convergence rates of the coupling scheme, which is again consistent with our theory. We point out that this observation may lead to significant improvement of the QM/MM coupling methods with {\it soft} matching condition, since the cost for training the virial stress is significantly less than training the MM models on the analogous order of derivatives.

\vskip 0.2cm

\subsubsection{Convergence results for the NRLTB Si model.}
We now move to testing our new QM/MM schemes when the QM reference model is an electronic structure model for Si. We chose NRLTB \cite{cohen94} as the reference model, which is a succesful tight-binding model for silicon. Our choice of Si as the material is due to the fact that it is a rich semi-conducting material for which we have also strong theoretical and numerical evidence for the localisation of its interatomic forces~\cite{chen18}.

\begin{figure}[!htb]
    \centering
    \includegraphics[scale=0.48]{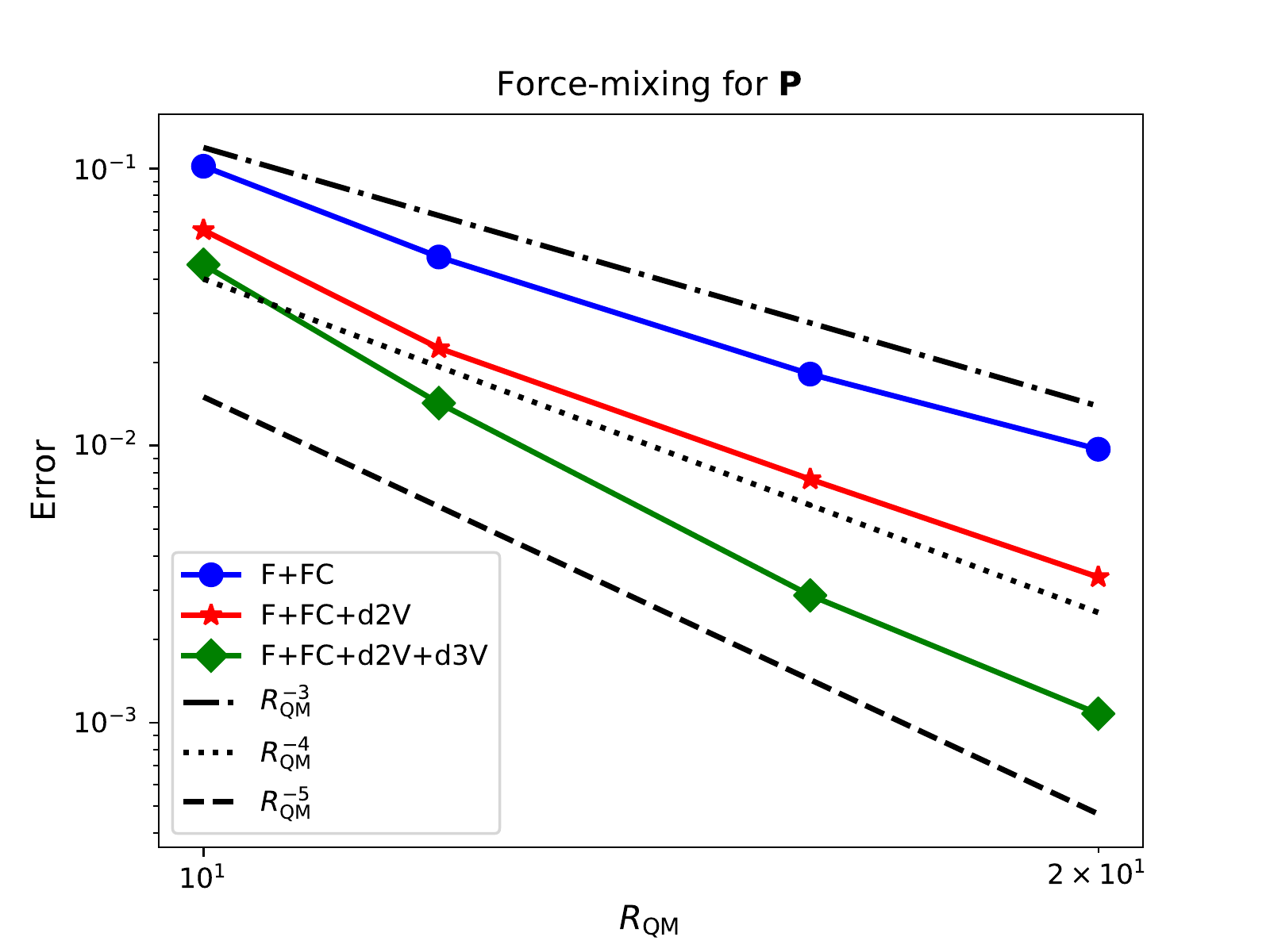}
    \includegraphics[scale=0.48]{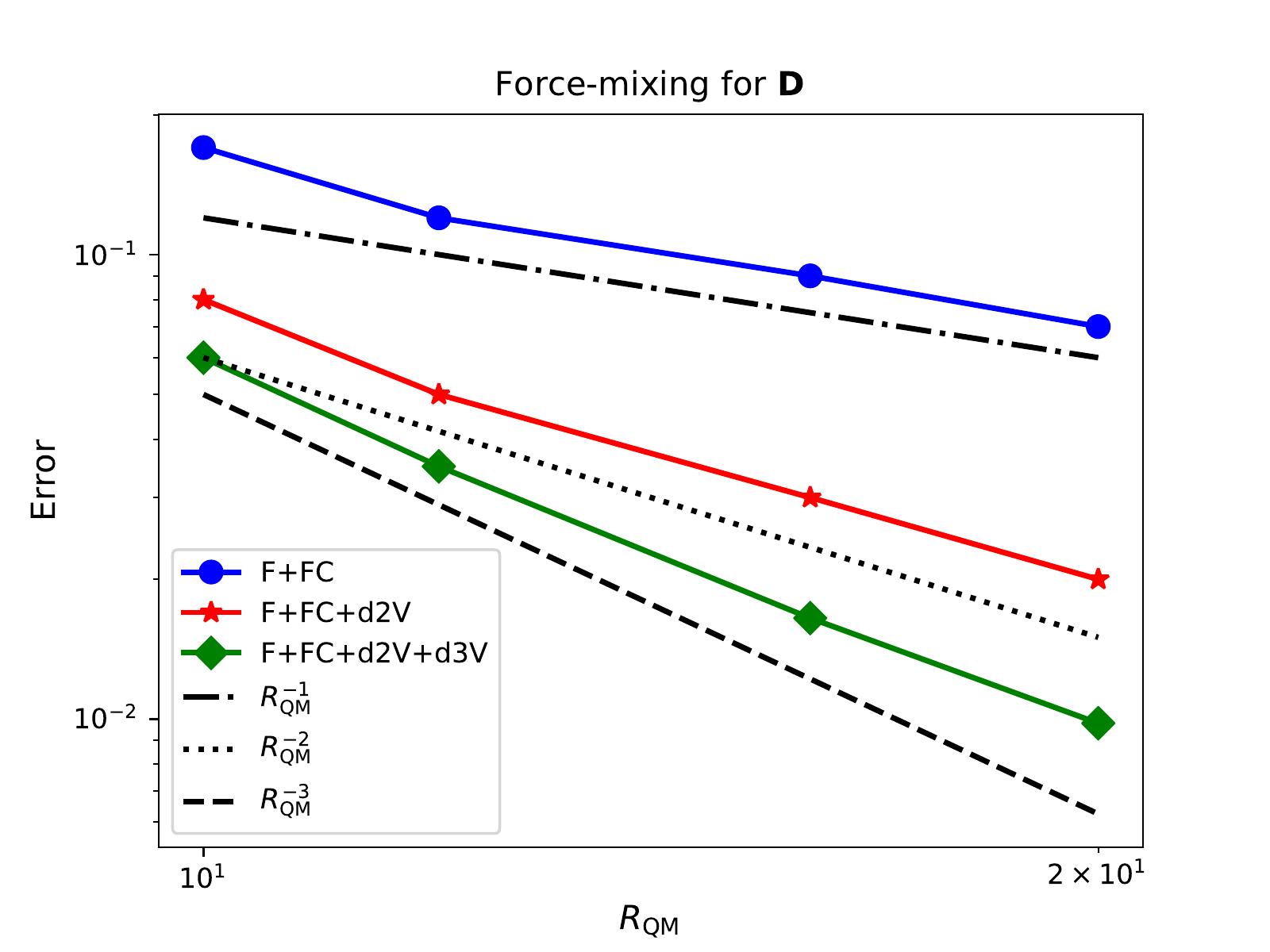}
    \caption{NRLTB Si reference model: Convergence rates of the force-mixing scheme for point defects (left) and dislocations (right).}
    \label{fig:force-mixing-TB}
\end{figure}

According to the discussions in the foregoing paragraph, the force-mixing schemes are far more efficient than the energy-mixing schemes. In the current section, we therefore only test the force-mixing schemes. The detailed fitting parameters, fitting accuracy in terms of RRMSE and fitting time of constructing MLIPs for force-mixing are given in Table \ref{table-params-TB}. 
\begin{table}
	\begin{center}
	\vskip0.5cm
	\begin{tabular}{|c|c|c|c|c|c|}
	\hline
		\diagbox{Schemes}{Accuracy (RRMSE)}{Measurements} & {$\ffit_1$} & {$\ffit_2$} & {$\vfit_2$} & {$\vfit_3$} & {T (s)}\\[0.5mm]
		\hline
		$K_{\rm F}=1$, $\# {\pmb O}=2400$, $\# {\pmb B}=383$  & 1.08\% & - & - & -  & 67.2 \\[0.5mm]
		\hline
		$K_{\rm F}=1~~\&$ virial,  $\# {\pmb O}=3129$, $\# {\pmb B}=383$  & 1.31\% & - & 0.48\% & - & 184.4
		\\[0.5mm]
	   \hline
	     $K_{\rm F}=1~~\&$ two virials, $\# {\pmb O}=9690$, $\# {\pmb B}=1119$ & 3.72\% & 5.03\% & - & 1.58\% & 2198.6
		\\[0.5mm]
		\hline
	\end{tabular}
	\end{center}
	\vskip 0.5cm
	\caption{NRLTB Si reference model: The fitting accuracy (RRMSE) on different measurements and fitting time for force-mixing (varying $K_{\rm F}$) coupling schemes.
	}
	\label{table-params-TB}
\end{table}
The $\RMM$ and $\Rbuf$ parameters are again chosen sufficiently large ($\RMM=80r_0$, $\Rbuf=8r_0$, where $r_0$ is the lattice constant of Si).
We then equilibrate our QM/MM models for (case \asP) a vacancy defect; and (case \asD) a (110)[100] straight edge-dislocation. 
The decay of approximate error $\|D\bar{u} - D\bar{u}^{\rm H}\|_{\ell^2_{\mathcal{N}}}$ of the QM/MM force-mixing schemes are shown in Figure \ref{fig:force-mixing-TB}.
We observe that the convergence rates again perfectly match our theoretical predictions from Theorems \ref{thm:forceM_V} and which are summarized in Table \ref{table-e-mix}. 
Moreover, it is worth mentioning that, the force-mixing scheme with {\it soft} matching conditions can achieve the third order convergence for straight edge-dislocations for realistic QM models without too much computational cost, while it is extremely difficult to reproduce this result by using the force-mixing scheme with {\it strict} matching conditions proposed in our previous paper \cite{chen17}. This observation makes it possible for us to simulate large scale QM/MM coupling methods for practical systems.


\section{Concluding remarks}
\label{sec:conclusion}
\setcounter{equation}{0}

In this paper, we propose and analyze new QM/MM coupling schemes for crystalline solids with embedded defects, where the MM models are machine-learned interatomic potentials (MLIPs).
We identify specific ``measurement'' to perform on the QM model which the MM model must match and provide a rigorous numerical analysis on the convergence in terms of these measurements and the sizes of QM and MM regions.
Our results not only explain the succes of existing QM/MM coupling methods (in particular force-mixing schemes) but also provide clear instructions on how to construct new interatomic potential (or force fields) specifically for the purpose of QM/MM coupling.

\appendix
\renewcommand\thesection{\appendixname~\Alph{section}}


\section{Far-field boundary predictor for dislocations}
\label{sec:appendixU0}
\renewcommand{\theequation}{A.\arabic{equation}}
\renewcommand{\thelemma}{A.\arabic{lemma}}
\setcounter{equation}{0}

\newcommand{\ulin}{u^{\rm lin}}
\newcommand{\burg}{{\sf b}}

We model dislocations by following the setting used in \cite{Ehrlacher16}. We consider a model for straight dislocations obtained by projecting a 3D crystal into 2D. Let $B \in \R^{3\times 3}$ be a nonsingular matrix. Given a Bravais lattice $B\Z^3$ with dislocation direction parallel to $e_3$ and Burgers vector $\burg=(\burg_1,0,\burg_3)$, we consider
displacements $W: B\Z^3 \rightarrow \R^3$ that are periodic in the direction of the dislocation direction of $e_3$. Thus, we choose a projected reference lattice $\L := A\Z^2 := \{(\ell_1, \ell_2) ~|~ \ell=(\ell_1, \ell_2, \ell_3) \in B\Z^3\}$. We also introduce the projection operator 
\begin{equation}\label{eq:P}
    P(\ell_1, \ell_2) = (\ell_1, \ell_2, \ell_3) \quad \text{for}~\ell \in B\Z^3.
\end{equation}
It can be readily checked that this projection is again a Bravais lattice.

We prescribe the far-field predictor $u_0$ as follows according to \cite{chen19, Ehrlacher16}. Let $\L\subset\R^2$, $\hat{x}\in\R^2$ be the position of the dislocation core and $\Gamma := \{x \in \R^2~|~x_2=\hat{x}_2,~x_1\geq\hat{x}_1\}$ be the ``branch cut'', with $\hat{x}$ chosen such that $\Gamma\cap\Lambda=\emptyset$.
We define the far-field predictor $u_0$ by
\begin{eqnarray}\label{predictor-u_0-dislocation}
u_0(x):=\ulin(\xi^{-1}(x)),
\end{eqnarray}
where $\ulin \in C^\infty(\R^2 \setminus \Gamma; \R^d)$ is the solution of continuum linear elasticity (CLE)
\begin{align}\label{CLE}
\nonumber
\mathbb{C}^{j\beta}_{i\alpha}\frac{\partial^2 u^{\rm lin}_i}{\partial x_{\alpha}\partial x_{\beta}} = 0 \qquad &\text{in} ~~ \R^2\setminus \Gamma,
\\
u^{\rm lin}(x+) - u^{\rm lin}(x-) = -\burg \qquad &\text{for} ~~  x\in \Gamma \setminus \{\hat{x}\},
\\
\nonumber
\nabla_{e_2}u^{\rm lin}(x+) - \nabla_{e_2}u^{\rm lin}(x-) = 0 \qquad &\text{for} ~~  x\in \Gamma \setminus \{\hat{x}\},
\end{align}
where the forth-order tensor $\mathbb{C}$ is the linearised Cauchy-Born tensor (derived from the potential $\Vhom$, see \cite[\S~7]{Ehrlacher16} for more detail),
\begin{eqnarray}
\xi(x)=x-\burg_{12}\frac{1}{2\pi}
\eta\left(\frac{|x-\hat{x}|}{\hat{r}}\right)
\arg(x-\hat{x}),
\end{eqnarray}
with $\arg(x)$ denoting the angle in $(0,2\pi)$ between $x$ and
$\burg_{12} = (\burg_1, \burg_2) = (\burg_1, 0)$, and
$\eta\in C^{\infty}(\R)$ with $\eta=0$ in $(-\infty,0]$ and $\eta=1$ in
$[1,\infty)$ which removes the singularity.

In order to model dislocations, the homogeneous site potential $V^{\rm h}$ must be invariant under lattice slip. Following \cite{Ehrlacher16}, we define the slip operator $S_0$ acting on the displacements $w: \Lambda \rightarrow \mathbb{R}^d$, by ($\mathsf{b}_{12}$ represents the projection of the Burger's vector to the $(x_1,x_2)$ plane)
\[
S_0w(x) := \begin{cases}
~~w(x) \qquad \qquad \qquad &x_2 > \hat{x}_2 \\
~~w(x-\mathsf{b}_{12}) - \mathsf{b} \qquad &x_2 < \hat{x}_2
\end{cases}.
\]

We may then formulate the slip invariance condition by defining a mapping $S$, where $S$ is an $\ell^2$-orthogonal operator with dual $S^{*}=S^{-1}$ by
\[
Su(\ell) := \begin{cases} ~~u(\ell) \qquad &\ell_2 > \hat{x}_2 \\
~~u(\ell-\mathsf{b}_{12}) \qquad &\ell_2 < \hat{x}_2
\end{cases}, \qquad 
S^{*}u(\ell) = \begin{cases} ~~u(\ell) \qquad &\ell_2 > \hat{x}_2 \\
~~u(\ell+\mathsf{b}_{12}) \qquad &\ell_2 < \hat{x}_2
\end{cases}.
\]
The slip invariance condition can now be expressed  as
\begin{equation}\label{slip}
V^{\rm h}\big(D(u_0 + u)(\ell)\big) = V^{\rm h}\big(S^{*}DS_0(u_0 + u)(\ell)\big), \qquad \forall~\ell \in \Lambda, u \in \UsH(\L),
\end{equation}
where $u_0$ is defined by \eqref{predictor-u_0-dislocation}.

In our analysis we require that applying the slip operator to the predictor map $u_0$ yields a smooth function in the half-space $\Omega_\Gamma = \{x_1 \geq \hat{x}_1 + \hat{r}+\mathsf{b}_{1}\}$.
It is therefore natural to define (likewise to \cite{Ehrlacher16}) the elastic strains
\begin{equation}\label{strain}
e(\ell) := \big(e_{\rho}(\ell)\big)_{\rho \in \L-\ell}, \qquad e_{\rho}(\ell) = \begin{cases} ~~S^{*}D_{\rho}S_0u_0(\ell) \qquad &\ell \in \Omega_{\Gamma} \\
~~D_{\rho}u_0(\ell) \qquad  &\ell \notin \Omega_\Gamma
\end{cases},
\end{equation}
and the analogous definition for corrector $u$
\begin{equation}\label{fancy_D}
Du(\ell) := \big(D_{\rho}u(\ell)\big)_{\rho \in \L-\ell}, \qquad D_{\rho}u(\ell) = \begin{cases} ~~S^{*}D_{\rho}Su(\ell) \qquad &\ell \in \Omega_{\Gamma} \\
~~D_{\rho}u(\ell) \qquad &\ell \notin \Omega_\Gamma
\end{cases}.
\end{equation}
Using this notation, the slip invariance condition \eqref{slip} may be written as, for $u\in\UsH(\L)$,
\begin{equation}\label{slip_strain}
V^{\rm h}\big(D(u_0 + u)(\ell)\big) = V^{\rm h}\big(e(\ell) + Du(\ell)\big).
\end{equation}

The following lemma, proven in \cite{chen19}, is a straightforward extension of \cite[Lemma 3.1]{Ehrlacher16}.
\begin{lemma}
If the predictor $u_0$ is defined by \eqref{predictor-u_0-dislocation} and $e(\ell)$ is given by \eqref{strain}, then there exists a constant $C$ such that 
\begin{equation}\label{eq:decay_u0}
    |e_{\sigma}(\ell)| \leq C |\sigma|\cdot |\ell|^{-1} \qquad {\rm and} \qquad |D_{\rho} e_{\sigma}(\ell)| \leq C |\rho| \cdot |\sigma|\cdot |\ell|^{-2}.
\end{equation}
\end{lemma}


\section{Proofs for the energy-mixing scheme}
\label{sec:appendixProofEmix}
\renewcommand{\theequation}{B.\arabic{equation}}\renewcommand{\thelemma}{B.\arabic{lemma}}
\renewcommand{\theremark}{B.\arabic{remark}}
\setcounter{equation}{0}

%
For $K\geq 2$, $\pmb{g}\in\big(\R^d\big)^{\LhomS}$, we define the $K$-th order Taylor's expansion $T_{K}\Vhom$ at the homogeneous state by
\begin{eqnarray}\label{eq:taylor_E}
T_{K}\Vhom\big({\bm g}\big) 
:= \Vhom({\bf 0}) + \sum_{j=1}^{K} \frac{1}{j!} \delta^j \Vhom({\bf 0})\left[{\bm g}^{\otimes j}\right].
\end{eqnarray}

In addition to the norm $\|D\cdot \|_{\ell^2_{\mathcal{N}}}$ introduced in \eqref{eq: nn norm}, one can construct a family of weighted norms that use weighted finite-difference stencils with infinite interaction range \cite{chen19}.  For a stencil $Du(\ell)$ and $\gamma>0$, we can define the (semi-)norms 
\begin{align}\label{eq:wnormdef}
\big|Du(\ell)\big|_{\wf_\gamma} := \bigg( \sum_{\rho\in\L-\ell}
\wf_\gamma(|\rho|) \big|D_\rho u(\ell)\big|^2 \bigg)^{1/2} \quad \text{and} \quad
\|Du\|_{\ell^2_{\wf_\gamma}} &:= \bigg( \sum_{\ell \in \L} |Du(\ell)|_{\wf_\gamma}^2 \bigg)^{1/2} 
\end{align}
with $\wf_\gamma(|\rho|):=e^{-2\gamma |\rho|}$. For convenience, we will use $\wf$ directly when there is no confusion in the context. We also have from \cite[Appendix A]{chen19} that for $\gamma>0$, there exist constants $c, C > 0$ such that
\begin{eqnarray}\label{eq:normest12}
 c\|Du\|_{\ell^2_{\mathcal{N}}} \leq \|Du\|_{\ell^2_{\wf}} \leq C \|Du\|_{\ell^2_{\mathcal{N}}} \qquad\forall~u\in\UsH(\L).
\end{eqnarray}


We are ready to give the proofs of the main results (Theorem \ref{theo:energyM} and Theorem \ref{thm:egyM_V}) for energy-mixing scheme and we essentially extend the proofs of \cite[Theorem 4.1]{chen17} to general cases.

\begin{proof}[{\bf Proof of Theorem \ref{theo:energyM}}]

Following the framework of the {\it a priori} error estimates in \cite{chen17, co2011, colz2012, colz2016}, we divide the proof into three steps in order to apply the inverse function theorem \cite[Lemma 2.2]{co2011}.

{\it 1. Quasi-best approximation:} We first construct $T^{\rm H} \bar{u}(\ell) \in \Admu^{\rm H}$ by
\begin{align}
T^{\rm H}\bar{u}(\ell) := \eta(\ell/\RMM)\big( \bar{u}(\ell) - a_{\RMM} \big),
\end{align}
where $\eta \in C^1(\R^d)$ is a cutoff function satisfying $\eta(x) = 1$ for $|x|\leq4/6$ and $\eta(x)=0$ for $|x|\geq5/6$ and $a_{\RMM}:=\bbint_{B_{5\RMM/6}\backslash B_{4\RMM/6}} I\bar{u}(x) \dx$ with $I\bar{u}$ defined by the piecewise affine interpolant of $\bar{u}$ with respect to the lattice \cite{colz2016, wang2020posteriori}. Then, for $\RMM$ sufficiently large, we have from the decay estimates \eqref{eq:ubar-decay} that 
\begin{align}
\label{proof-TH-diff}
~&\|D\bar{u}-DT^{\rm H}\bar{u}\|_{\ell^2_{\wf}} \leq C \|D\bar{u}\|_{\ell^2_{\wf}(\L\backslash B_{\RMM/2})}\qquad \text{and} \\
\label{proof-TH}
&|DT^{\rm H}\bar{u}(\ell)|_{\wf} \leq C(1+|\ell|)^{-d}\log^t(2+|\ell|) \quad \forall~\ell \in \L.
\end{align}
Using the norm equivalence \eqref{eq:normest12}, for $\RMM$ sufficiently large, we have
        \begin{align}
                \label{proof-4-1-2}
 \|\delta\E(\bar{u})-\delta\E(T^{\rm H}\bar{u})\|
                &\leq L_1\|D\bar{u}-DT^{\rm H}\bar{u}\|_{\ell^2_{\mathcal{N}}}
                \leq CL_1\|D\bar{u}\|_{\ell^2_{\mathcal{N}}(\L\backslash B_{\RMM/2})}
                \qquad \text{and} \\
                \label{proof-4-1-3}
                \|\delta^2\E(\bar{u})-\delta^2\E(T^{\rm H}\bar{u})\|
                &\leq L_2\|D\bar{u}-DT^{\rm H}\bar{u}\|_{\ell^2_{\mathcal{N}}}
                \leq CL_2\|D\bar{u}\|_{\ell^2_{\mathcal{N}}(\L\backslash B_{\RMM/2})},
        \end{align}
where $L_1$ and $L_2$ are uniform Lipschitz constants of $\delta\E$ and $\delta^2\E$ respectively since $\E$ is $(\n-1)$-times continuously differentiable with respect to the $\|D\cdot\|_{\ell^2_{\mathcal{N}}}$ norm.

{\it 2. Stability:} 
We first observe that the ghost force correction term in \eqref{eq:hybrid_energy_BGFC} is a linear term with respect to $u$, hence its second variation with respect to $v\in\Admu^{\rm H}$ vanishes, that is, 
\begin{eqnarray}\label{eq:stab_eq}
\big \<\delta^2 \E^{\rm GFC}(T^{\rm H}\bar{u})v, v\big \> = \big \<\delta^2 \EH(T^{\rm H}\bar{u})v, v\big \> \quad \forall~v \in \Admu^{\rm H}.
\end{eqnarray}
Hence, for any $v \in \Admu^{\rm H}$, we consider the stability of
\begin{eqnarray}\label{eq:stab_hybrid}
\nonumber
&&\<\delta^2 \EH(T^{\rm H}\bar{u})v, v\> 
\\[1ex]
\nonumber 
&=& \big \<\delta^2 \E^{\rm T}(T^{\rm H}\bar{u})v, v\big \> + \Big( \big \<\delta^2 \EH(T^{\rm H}\bar{u})v, v\big \> - \big \<\delta^2 \E^{\rm T}(T^{\rm H}\bar{u})v, v\big \> \Big)
\\[1ex]
&=:& S_1 + S_2,
\end{eqnarray}
where $\E^{\rm T}$ is the hybrid energy functional 
\begin{align}\label{eq:hybrid_taylor_energy}
\quad \E^{\rm T}(u) 
:= \sum_{\ell\in \LQM}& 
\Big( \VQM_{\ell}\big(Du_0(\ell) + D u(\ell)\big) - \VQM_{\ell}\big(Du_0(\ell)\big) \Big)
\nonumber \\
+& \sum_{\ell\in \LMM\cup\LFF} 
\Big( T_{K} V^{\rm BUF}_{\#}\big(Du_0(\ell) + Du(\ell)\big) - T_{K} V^{\rm BUF}_{\#}\big(Du_0(\ell)\big) \Big),
\end{align}
where $V^{\rm BUF}_{\#}: (\R^d)^{\mathcal{R}} \rightarrow \R$ with $\mathcal{R}=B_{\Rbuf} \cap (\L \setminus 0)$
satisfies $V^{\rm BUF}_{\#}\big(D_{\mathcal{R}}u(\ell)\big) = V_{\ell}^{B_{\Rbuf}(\ell)}\big(Du(\ell)\big)$ and $T_K V^{\rm BUF}_{\#}$ is $K$-order Taylor's expansion of $V^{\rm BUF}_{\#}$. The term $S_1$ can be estimated by using the results of \cite[Theorem 4.1]{chen17}
\begin{eqnarray}\label{eq:stab_taylor}
 S_1 := \big\<\delta^2 \E^{\rm T}(T^{\rm H}\bar{u})v, v\big \> \geq c_{\rm E} \|Dv \|^2_{\ell^2_{\mathcal{N}}}
\end{eqnarray}
with some constants $c_{\rm E}>0$ that are independent of model parameters.

To estimate $S_2$, for simplicity of notation, we denote $\uh(\ell):= u_0(\ell) + I^{\rm h}\THu(\ell)$ and $\vh(\ell):=I^{\rm h}v(\ell)$ for $v\in\Admu^{\rm H}$ with $I^{\rm h}$ defined by \eqref{eq:eqI}, we then split it into three parts
\begin{eqnarray}\label{eq:upper_stab}
\nonumber
S_2 &=& \big \<\delta^2 \EH(T^{\rm H}\bar{u})v, v\big \> - \big \<\delta^2 \E^{\rm T}(T^{\rm H}\bar{u})v, v\big \>
\\[1ex]
\nonumber
&=& \sum_{\ell \in \LMM \cup \LFF} \Big\<\Big( \delta^2 \VMM\big(D\uh(\ell)\big) - \delta^2 T_{K}\VMM\big(D\uh(\ell)\big) \Big) D\vh(\ell), D\vh(\ell) \Big\>  
\\
\nonumber
&& + \sum_{\ell \in \LMM \cup \LFF} \Big\<\Big(\delta^2 T_{K}\VMM\big(D\uh(\ell)\big) - \delta^2 T_{K} V^{\rm h}\big(D\uh(\ell)\big) \Big) D\vh(\ell), D\vh(\ell) \Big\>  
\\
\nonumber
&& + \sum_{\ell \in \LMM \cup \LFF} \Big( \big\<\delta^2 T_{K} V^{\rm h}\big(D\uh(\ell)\big) D\vh(\ell), D\vh(\ell) \big\>  
\\
\nonumber
&& \qquad \qquad \qquad - \big\< \delta^2 T_{K} V^{\rm BUF}_{\#}\big(Du_0(\ell) + D\THu(\ell)\big) Dv(\ell), Dv(\ell) \big\>  \Big) 
\\[1ex]
&=:& S_{21} + S_{22} + S_{23}.
\end{eqnarray}
Using Taylor's expansion of $\VMM$ at the homogeneous state, we have 
        \begin{align}\label{eq:s11}
                |S_{21}| =& \sum_{\ell\in\LMM\cup\LFF}
                \left\< \Big( \delta^2 \VMM\big(D \uh(\ell)\big)
                - \delta^2 T_K \VMM\big(D \uh(\ell)\big) \Big)
                 D \vh(\ell) , D \vh(\ell) \right\>
                \nonumber \\
                \leq&~C\sum_{\ell\in\LMM\cup\LFF}
                |D \uh(\ell)|_{\wf}^{K-1} |D \vh(\ell)|^2_{\wf}
                \leq C\|D \uh\|^{K-1}_{\ell^{\infty}_{\wf}(\LMM\cup\LFF)}
                \|Dv\|^2_{\ell^2_{\wf}}.
        \end{align}
For $S_{22}$, we have
\begin{eqnarray}\label{eq:s12}
\nonumber
|S_{22}| &=& \sum^{K}_{j = 2}\sum_{\ell \in \LMM \cup \LFF} \Big\<\Big(\delta^{j} \VMM(\p0) - \delta^{j} \Vhom(\p0) \Big)\big(D\uh(\ell)\big)^{j-2} D\vh(\ell), D\vh(\ell) \Big\>
\\ \nonumber
&\leq& \efit_2 \sum_{\ell \in \LMM \cup \LFF} |D \vh(\ell)|^2_{\wf} + \sum^{K}_{j=3} \efit_{j} \sum_{\ell \in \LMM \cup \LFF} |D \uh(\ell)|_{\wf}^{j-2} |D \vh(\ell)|^2_{\wf}
\\
&\leq& \Big( \efit_2 + \sum^{K}_{j=3} \efit_{j} \|D \uh\|^{j-2}_{\ell^{\infty}_{\wf}(\LMM\cup\LFF)} \Big) \|Dv\|^2_{\ell^2_{\wf}}.
\end{eqnarray}  
%
Note that the locality assumption in \assERL~implies that there exists a constant $\eta>0$ such that
\begin{eqnarray}\label{eq:s13}
|S_{23}| \leq C e^{-\eta \Rbuf} \|Dv \|^2_{\ell^2_{\wf}}.
\end{eqnarray}

Combing from \eqref{eq:stab_eq} to \eqref{eq:s13} and using the decay estimates \eqref{eq:ubar-decay} and the norm equivalence \eqref{eq:normest12}, we have that, for sufficiently large $\RQM$ and sufficiently small $\efit_2$, 
\begin{eqnarray}\label{eq:stab}
\big \<\delta^2 \E^{\rm GFC}(T^{\rm H}\bar{u})v, v\big \>  \geq \frac{c_{\rm E}}{2} \|Dv \|^2_{\ell^2_{\mathcal{N}}}.
\end{eqnarray}

{\it 3. Consistency:} For any $v\in \Admu^{\rm H}$, we have
\begin{eqnarray}\label{eq:proof-c}
\nonumber
&& \big \<\delta \E^{\rm GFC}(\THu), v\big\>
\\[1ex]
\nonumber
&=& \big \<\delta \E^{\rm GFC}(\THu) - \delta \EH(\THu), v\big\> + \big \<\delta \EH(\THu) - \delta \E(\THu), v\big\> 
\\
\nonumber
&& \qquad + \big \<\delta \E(\THu) - \delta \E(\bar{u}), v\big\>
\\[1ex]
\nonumber
&=& 
\sum_{\ell\in\LQM} \big\< \delta \VQM_{\ell}\big(Du_0(\ell) + DT^{\rm H}\bar{u}(\ell)\big)
                - \delta V_{\ell}\big(Du_0(\ell) + DT^{\rm H}\bar{u}(\ell)\big), D v(\ell) \big\>
\\
\nonumber
&& \qquad + \sum_{\ell\in\LMM\cup\LFF} \Big( \big\< \delta \Vhom\big(D \uh(\ell)\big), D \vh(\ell) \big\> - \big\< \delta V_{\ell}\big(D u_0(\ell) + D T^{\rm H}\bar{u}(\ell)\big), D v(\ell) \big\> \Big)
\\[1ex]
\nonumber
&& + \sum_{\ell\in\LMM\cup\LFF} \big\< \delta \VMM\big(D \uh(\ell)\big)
                - \delta T_K \VMM\big(D \uh(\ell)\big), D \vh(\ell) \big\>
\\
\nonumber
&& \qquad +
\sum_{\ell\in\LMM\cup\LFF}
                \big\< \delta T_K\Vhom\big(D\uh(\ell)\big)
                - \delta  \Vhom\big(D \uh(\ell)\big), D \vh(\ell) \big\> 
\\[1ex]
\nonumber
&& + \quad \big \<\delta \E(\THu) - \delta \E(\bar{u}), v\big\>
\\[1ex]
\nonumber
&& + \sum_{\ell\in\LMM\cup\LFF}
                \big\< \delta T_K\VMM\big(D \uh(\ell)\big)
                - \delta T_K \Vhom\big(D \uh(\ell)\big), D \vh(\ell) \big\> - \big \< \delta \EH(\p0), \beta v \big \>
\\[1ex]
&=:& T_1 + T_2 + T_3 + T_4.
\end{eqnarray}
Following from the locality assumptions in \assERL~and the definition of the interpolation operator $I^{\rm h}$ by \eqref{eq:eqI}, there exits a constant $\kappa_1 > 0$ such that 
\begin{eqnarray}\label{proof-c-T1}
                |T_1| 
                &\leq& C e^{-\kappa_1 \Rbuf}
                \|Dv\|_{\ell^2_{\wf}}.
\end{eqnarray}
To estimate $T_2$, we have from \eqref{eq:taylor_E} and \eqref{eq:eqI}
\begin{eqnarray}\label{proof-c-T2}
\nonumber
|T_2| &=& \sum_{\ell\in\LMM\cup\LFF}
                \big\< \delta \VMM\big(D \uh(\ell)\big)
                - \delta T_K\VMM \big(D \uh(\ell)\big) ,
                D \vh(\ell) \big\>
                \\ \nonumber 
 \qquad \qquad && + \sum_{\ell\in\LMM\cup\LFF}
                \big\< \delta T_K\Vhom\big(D\uh(\ell)\big)
                - \delta  \Vhom\big(D \uh(\ell)\big), D \vh(\ell) \big\>
                \nonumber \\[1ex]
                &\leq& C\sum_{\ell\in\LMM\cup\LFF}
                |D \uh(\ell)|_{\wf}^{K} |D \vh(\ell)|_{\wf} \leq C\|D \uh\|^{K}_{\ell^{2K}_{\wf}(\LMM\cup\LFF)}
                \|Dv\|_{\ell^2_{\wf}}.
        \end{eqnarray}
For $T_3$, after a direct calculation, we have from \eqref{proof-4-1-2}
        \begin{eqnarray}\label{proof-c-T4}
                |T_3| \leq CL_1
                \|D\bar{u}\|_{\ell^2_{\wf
                }(\L\backslash B_{\RMM/2})}
                \|Dv\|_{\ell^2_{\wf}}.
        \end{eqnarray}
To estimate $T_4$, we first split it into two parts
\begin{eqnarray}\label{eq:proof-c-T5}
\nonumber
T_4 &=& \sum_{\ell\in\LMM\cup\LFF}
                \big\< \delta T_K\VMM\big(D\uh(\ell)\big)
                - \delta T_K \Vhom\big(D\uh(\ell)\big), D \vh(\ell) \big\> - \big \< \delta \EH(\p0), \beta v \big \> 
\\[1ex] \nonumber
&=& \sum_{\ell\in\LMM\cup\LFF} \big\< \delta \VMM(\p0)
                - \delta \Vhom(\p0), D \vh(\ell) \big\> - \big \< \delta \EH(\p0), \beta v \big \> 
\\ \nonumber
&& + \sum_{j=2}^K\sum_{\ell \in \LMM \cup \LFF} \Big\< \Big( \delta^{j} \VMM (\p0) - \delta^{j} \Vhom (\p0) \Big) \big(D\uh(\ell)\big)^{j-1} , D\vh(\ell) \Big\>
\\[1ex]
&=:& T_{41} + T_{42}.
\end{eqnarray}
To estimate $T_{41}$, for simplicity of notations, let $v_{\beta}(\ell):=v(\ell)\beta(\ell)$, $\Omega^{{\rm h}, c}_{\rm QM}:= \{\ell \in \Lhom ~|~ |\ell| \leq c\RQM \}$ with the constant $0 \leq c \leq 1$. It is straightforward to write $\Omega^{{\rm h}, 1/2}_{\rm QM} \subset \Omega^{{\rm h}, 1}_{\rm QM}$, then we have
\begin{eqnarray}\label{eq:T42_attep}
\nonumber
T_{41} 
&=& \sum_{\ell \in \Omega^{\rm h, 1}_{\rm QM}} \F^{\rm h}_{\ell}(\p0)\vh(\ell) - \sum_{\ell \in \L_{\rm QM}} \F^{\rm QM}_{\ell}(\p0)v_{\beta}(\ell)
\\ \nonumber
&& + \sum_{\ell \in \LMM \cup \LFF} \sum_{\substack{\rho\in\LhomS,\\ \ell+\rho \in \Omega_{\rm QM}^{{\rm h},1}}} \Big( \VQM_{\ell, \rho}(\p0)v(\ell+\rho) - \Vhom_{, \rho}(\p0)\vh(\ell+\rho) \Big)
\\ \nonumber
&&  + \sum_{\ell \in \LMM \cup \LFF}\sum_{\substack{\rho\in\LhomS,\\ \ell+\rho \in \Omega_{\rm QM}^{{\rm h}, 1/2}}} \VMM_{, \rho}(\p0) \vh(\ell+\rho)
\\[1ex]
&=:& T^{\rm (a)}_{41} + T^{\rm (b)}_{41}.
\end{eqnarray}  
The term $T^{\rm (a)}_{41}$ can be divided into two parts according to the blending function $\beta$, and noticing $\F^{\rm h}_{\ell}(\p0)=0$ for every $\ell \in \Lhom$, after a direct calculation, then we have
\begin{eqnarray}\label{ea:T521}
\nonumber
T^{\rm (a)}_{41} &=& \sum_{|\ell| \leq \RQM/2} \F^{\rm h}_{\ell}(\p0)\vh(\ell) + \sum_{\RQM/2 \leq |\ell| \leq \RQM} \big( \F^{\rm h}_{\ell}(\p0) - \F^{\rm QM}_{\ell}(\p0) \big) v(\ell)
\\[1ex] 
&\leq& Ce^{-\kappa_2 \Rbuf} \|Dv\|_{\ell^2_{\wf}},
\end{eqnarray}
where the last inequality follows from the locality assumption in \assERL.
The term $T^{\rm (b)}_{41}$ can be estimated similarly 
\begin{eqnarray}\label{eq:T41b}
|T^{\rm (b)}_{41}| \leq C e^{-\kappa_3 \Rbuf} \|Dv\|_{\ell^2_{\wf}}.
\end{eqnarray}
We then estimate $T_{42}$ by
\begin{eqnarray}\label{eq:T51}
\nonumber
|T_{42}| &=& \sum_{j=2}^K \sum_{\ell \in \LMM \cup \LFF} \Big\< \Big( \delta^{j} \VMM (\p0) - \delta^{j} \Vhom (\p0) \Big) \big(D\uh(\ell)\big)^{j-1} , D\vh(\ell) \Big\> 
\\[1ex] \nonumber
&\leq& \sum_{j=2}^K \efit_{j} \sum_{\ell \in \LMM \cup \LFF}  \big|D \uh(\ell)\big|_{\wf}^{j-1} \big|D \vh(\ell)\big|_{\wf} \leq \Big( \sum^{K}_{j=2} \efit_{j} \|D \uh\|^{j-1}_{\ell^{2j-2}_{\wf}(\LMM\cup\LFF)} \Big) \|Dv\|_{\ell^2_{\wf}}.
\end{eqnarray}  
Taking into accounts from \eqref{eq:proof-c} to \eqref{eq:T41b}, and the norm equivalence \eqref{eq:normest12}, let $\kappa:=\max\{\kappa_1, \kappa_2, \kappa_3\}$, we have
	\begin{align}\label{e-mix-consistency}
		\big\< \delta\E^{\rm GFC}(T^{\rm H}\bar{u}), v \big\> \leq~& C
		\Big( \sum^{K}_{j=2} \efit_{j} \|D\uh\|^{j-1}_{\ell^{2j-2}_{\wf}(\LMM\cup\LFF)} + \|D\uh\|^{K}_{\ell^{2K}_{\wf}(\LMM\cup\LFF)}
		\nonumber \\
		& \qquad  + \|D\bar{u}\|_{\ell^{2}_\wf(\L\backslash B_{\RMM/2})} + e^{-\kappa \Rbuf}
		\Big) \cdot \|Dv\|_{\ell^2_{\mathcal{N}}}.
	\end{align}

For case {\bf (P)}, we obtain from $u_0={\p0}$ and $|D \bar{u}(\ell)|_\wf \leq C (1+|\ell|)^{-d}$
	\begin{align}\label{e-mix-consistency-P}
		\left| \big\< \delta\E^{\rm GFC}(T^{\rm H}\bar{u}) , v \big\> \right| \leq C
		\Big( \sum^{K}_{j=2} \efit_j  R^{-(2j-3)d/2}_{\rm QM} + \RQM^{-(2K-1)d/2} 
		+ \RMM^{-d/2} + e^{-\kappa\Rbuf} \Big) \cdot
		\|Dv\|_{\ell^2_{\mathcal{N}}}.
	\end{align}

For case {\bf (D)}, we obtain from \eqref{eq:decay_u0}, $|Du_0(\ell)|_{\wf}\leq C|\ell|^{-1}$ and $|D\bar{u}(\ell)|_{\wf}\leq C(1+|\ell|)^{-2}\log(2+|\ell|)$	\begin{align}\label{e-mix-consistency-D}
		\left| \big\< \delta\E^{\rm GFC}(T^{\rm H}\bar{u}) , v \big\> \right| \leq C
		\Big( \sum^{K}_{j=2} \efit_j  R^{-j+2}_{\rm QM} + \RQM^{-K+1} 
		+ \RMM^{-1}\log\RMM + e^{-\kappa\Rbuf}\Big) \cdot
		\|Dv\|_{\ell^2_{\mathcal{N}}}.
	\end{align}

Combining the stability \eqref{eq:stab} with consistency \eqref{e-mix-consistency-P} and \eqref{e-mix-consistency-D}, we can apply the inverse function theorem \cite[Lemma 2.2]{co2011} with $K=K_{\rm E}$ and use \eqref{proof-TH-diff} to obtain the stated results.
\end{proof}

\begin{remark}\label{re:gfc}
If the ghost force correction \eqref{eq:hybrid_energy_BGFC} is not applied, that is, $\beta \equiv 0$, then the estimate of $T_{41}$ in \eqref{eq:T42_attep} becomes to
\begin{eqnarray}\label{eq:T52_nogfc}
\nonumber
|T_{41}| 
&\leq& \efit_1 \sum_{\ell \in \LMM \cup \LFF} \big|D \vh(\ell)\big|_{\wf} \leq C \efit_1 \cdot N_{\rm MM}^{1/2} \cdot \|Dv\|_{\ell^2_{\wf}},
\end{eqnarray}  
which is non-vanishing as $\efit_1$ is not exact zero (slight mismatch at the QM/MM interface). It gives rise to the so-called ``ghost-force" inconsistency.  
Hence the ghost force correction term $\<\delta \EH(\p0), \beta v\>$ is required to reduce the effect of ``ghost-force". 
\end{remark}

We then give the proof of Theorem \ref{thm:egyM_V}, which is similar to that of Theorem \ref{theo:energyM}. Due to adding the soft matching conditions for virials in Theorem \ref{thm:egyM_V}, the only difference lies in the consistency part.
\begin{proof}[{\bf Proof of Theorem \ref{thm:egyM_V}}]
The quasi-best approximation and the stability analysis are the same as that in the proof of Theorem \ref{theo:energyM}. The main difference is the estimate of $T_2$ in \eqref{eq:proof-c}. $T_2$ can be further splitted into two parts
\begin{eqnarray}\label{proof-c-T2-v}
\nonumber
T_2 &=& \sum_{\ell\in\LMM\cup\LFF}
                \big\< \delta \VMM\big(D \uh(\ell)\big)
                - \delta T_K\VMM \big(D \uh(\ell)\big) ,
                D \vh(\ell) \big\>
                \\ \nonumber 
 \qquad \qquad && + \sum_{\ell\in\LMM\cup\LFF}
                \big\< \delta T_K\Vhom\big(D\uh(\ell)\big)
                - \delta  \Vhom\big(D \uh(\ell)\big), D \vh(\ell) \big\>
                \nonumber \\[1ex]
&=:& \sum_{\ell\in\LMM\cup\LFF} \big\<\mathscr{G}_1\big(D \uh(\ell)\big), D \vh(\ell) \big\> + \big\<\mathscr{G}_2\big(D \uh(\ell)\big), D \vh(\ell) \big\>,
\end{eqnarray}
where $\mathscr{G}_1:=\delta\VMM-\delta T_{K+1}\VMM + \delta T_{K+1} \Vhom - \delta \Vhom$ and $\mathscr{G}_2:=\delta T_{K+1}\VMM-\delta T_{K}\VMM + \delta T_{K} \Vhom - \delta T_{K+1} \Vhom$.
The first term of \eqref{proof-c-T2-v} can be similarly estimated by using the techniques in \eqref{proof-c-T2}
\begin{equation}\label{eq:T2-v-1}
\Big| \sum_{\ell\in\LMM\cup\LFF} \big\<\mathscr{G}_1\big(D \uh(\ell)\big), D \vh(\ell) \big\> \Big| \leq C \|D \uh\|^{K+1}_{\ell^{2K+2}_{\wf}(\LMM\cup\LFF)}
                \|Dv\|_{\ell^2_{\wf}}.
\end{equation}
For the second term, let $\Omega^{\rm MF}:= \{x \in \R^d, |x| \geq \RQM\}$, and expanding the strain $D_{\rho}u(\ell)$ at $\ell$, then we have
\begin{eqnarray}\label{eq:T2-v-2}
\nonumber
&& \Big| \sum_{\ell\in\LMM\cup\LFF} \big\<\mathscr{G}_2\big(D \uh(\ell)\big), D \vh(\ell) \big\> \Big|
\nonumber \\
&\leq& \sum_{\ell\in\LMM\cup\LFF}\sum_{\pmb{\rho} \in (\LhomS)^{K+1}}\left( \Vhom_{, \pmb{\rho}}(\p0) - \VMM_{, \pmb{\rho}}(\p0) \right) \Big( \prod^{K}_{i=1} \nabla \uh(\ell) \cdot \rho_i \Big) \big( \nabla \vh(\ell) \cdot \rho_{K+1} \big)
\nonumber \\
\qquad \qquad && + C\sum_{\ell\in\LMM\cup\LFF}
                \big(\nabla \uh(\ell)\big)^{K-1} \nabla^2\uh(\ell)\nabla  \vh(\ell)
\nonumber \\[1ex]
&\leq& C \big( \vfit_K \|\nabla \uh\|^{K}_{L^{2K}(\Omega^{\rm MF})} + \|\nabla \uh\|^{K-1}_{L^{2K-2}(\Omega^{\rm MF})}\|\nabla^2 \uh\|_{L^2(\Omega^{\rm MF})}\big) \cdot \|Dv\|_{\ell^2_{\wf}},
\end{eqnarray}
where the last inequality follows from the Cauchy-Schwarz inequality and the norm equivalence between $\|Dv\|_{\ell^2_{\wf}}$ and $\|\nabla v\|_{L^2}$.

Hence, substituting the estimate for $T_2$ in \eqref{eq:proof-c} with the new estimates (combining \eqref{proof-c-T2-v}, \eqref{eq:T2-v-1} and \eqref{eq:T2-v-2}), we have the following consistency results
	\begin{align}\label{e-mix-consistency-v}
		\big\< \delta\E^{\rm GFC}(T^{\rm H}\bar{u}), v \big\> \leq~& C
		\Big( \sum^{K}_{j=2} \efit_{j} \|D\uh\|^{j-1}_{\ell^{2j-2}_{\wf}(\LMM\cup\LFF)} + \vfit_K \|\nabla \uh\|^{K}_{L^{2K}(\Omega^{\rm MF})}
		\nonumber \\
		&+ \|\nabla \uh\|^{K-1}_{L^{2K-2}(\Omega^{\rm MF})}\|\nabla^2 \uh\|_{L^2(\Omega^{\rm MF})} + \|D\bar{u}\|_{\ell^{2}_\wf(\L\backslash B_{\RMM/2})} + e^{-\kappa \Rbuf}
		\Big) \cdot \|Dv\|_{\ell^2_{\mathcal{N}}}.
	\end{align}
	
Combining the decay estimates \eqref{eq:ubar-decay}, the stability \eqref{eq:stab} and \eqref{proof-TH-diff}, we can obtain the stated results.
\end{proof}


\section{Proofs for the force-mixing scheme}
\label{sec:appendixProofFmix}
\renewcommand{\theequation}{C.\arabic{equation}}
\renewcommand{\thelemma}{C.\arabic{lemma}}
\setcounter{equation}{0}

For $v:\L\rightarrow\R$, we will use the notation $\big\<\FH(u),v\big\>:=\sum_{\ell\in\LQM\cup\LMM}\FH_{\ell}(u)\cdot v(\ell)$ in our analysis.
For $K \geq 1$ and $w:\L \rightarrow \R^d$, we introduce the $K$-th order Taylor's expansion $T_K\F^{\rm h}$ at the homogeneous state by 
\begin{eqnarray}\label{eq:taylor_F}
T_K\F^{\rm h}\big(w\big) 
:= \sum_{j=1}^K \frac{1}{j!} \delta^j \F^{\rm h}({\bf 0})\left[w^{\otimes j}\right].
\end{eqnarray}

The following two Lemmas are useful in our analysis in that they sometimes allow us to avoid stress-strain representations of residual forces that we need to estimate.  We include them here for the sake of completeness and please refer to \cite{chen17} for more details.
\begin{lemma}\label{le:v}
For any $v \in \Admu^{\rm H}$, there exists $C>0$ such that
\begin{align}
    \|v\|_{\ell^{\infty}} &\leq C \|Dv\|_{\ell^2_{\wf}}(1+\log\RMM) \quad \textrm{\rm if}~d=2 ~\textrm{\rm and} 
    \\
    \|v\|_{\ell^6} &\leq C \|Dv\|_{\ell^2_{\wf}} \qquad  \qquad  \qquad \quad \text{\rm if}~d=3.
\end{align}
\end{lemma}
\begin{lemma}\label{le:fv}
Let $0<L<R$ and $v:\Lambda \rightarrow \R$ satisfy $v(\ell)=0 ~ \forall ~ |\ell| \geq R$. If $f: \Lambda \rightarrow \R$ satisfies $|f(\ell)| \leq c|\ell|^{-p}$ with $p\geq 3$, then there exists a constant $C$ such that 
\begin{equation}\label{le:div}
    \sum_{L\leq|\ell|\leq R} f(\ell) \cdot v(\ell) \leq C \Big(1+\log\Big(\frac{R}{L}\Big)\Big) \cdot L^{-p+1+d/2} \cdot \|Dv\|_{\ell^2_{\wf}} \qquad \forall~v \in \Admu^{\rm H}. 
\end{equation}
\end{lemma}

We are in positions to give the proofs of the main results (Theorem \ref{theo:forceM} and Theorem \ref{thm:forceM_V}) for force-mixing scheme and most of the proofs essentially follow from \cite[Proof of Theorem 5.1]{chen17}.

\begin{proof}[{\bf Proof of Theorem \ref{theo:forceM}}]

Similar to the proofs for energy-mixing, we divide the proof into three steps in order to apply the inverse function theorem \cite[Lemma 2.2]{co2011}.

{\it 1. Quasi-best approximation:} We take the approximation $T^{\rm H}\bar{u}\in\Admu^{\rm H}$ constructed in \ref{sec:appendixProofEmix}, so that the properties from \eqref{proof-TH-diff} to \eqref{proof-4-1-3} are satisfied.

{\it 2. Stability:} For any $v \in \Admu^{\rm H}$, we have
\begin{eqnarray}\label{eq:stab_hybrid_frc}
\nonumber
&&\<\delta \FH(T^{\rm H}\bar{u})v, v\>
\\[1ex]
\nonumber 
&=& \big \< \delta \F^{\rm T}(T^{\rm H}\bar{u})v, v \big\> + \big \< \big (\delta \FH(T^{\rm H}\bar{u}) - \delta \F^{\rm T}(T^{\rm H}\bar{u})\big)v, v\big \> 
\\[1ex]
&=:& G_1 + G_2,
\end{eqnarray}
where $\F^{\rm T}$ is the QM/MM hybrid force 
\begin{align}\label{F-T}
	\F_{\ell}^{\rm T}(u) =
	\left\{ \begin{array}{ll}
		\F^{\LQM \cup \Lbuf}_{\ell}(u) \qquad & \ell\in\LQM,
		\\[1.5ex]
		T_K \F^{\rm BUF}_{\#}\left( \big( u_0(\cdot-\ell) +  u(\cdot-\ell)\big)
		\big|_{B_{\Rbuf}} \right) \qquad & \ell\in\LMM,  
	\end{array} \right. 
\end{align}
where $\F^{\rm BUF}_{\#}: (\R^d)^{\mathcal{R}} \rightarrow \R$ with $\mathcal{R}=B_{\Rbuf} \cap \Lhom$
satisfying
\begin{eqnarray*}
\F^{\rm BUF}_{\#}\big(u_0(\cdot-\ell) + u(\cdot - \ell)\big|_{B_{\Rbuf}}\big) = \F_{\ell}^{B_{\Rbuf}(\ell)}\big(u\big).
\end{eqnarray*} 

The term $G_1$ can be estimated by using \cite[Theorem 5.1]{chen17} directly
\begin{align}\label{eq:stab_taylor_frc}
G_1 := \big \<\delta \F^{\rm T}(T^{\rm H}\bar{u})v, v\big \> \geq c_{\rm F} \|Dv \|^2_{\ell^2_{\wf}},
\end{align}
with some constants $c_{\rm F}>0$ that are independent of model parameters.

To estimate $G_2$, for simplicity of notations, let $\uh(\cdot-\ell):=u_0(\cdot-\ell) + I^{\rm h}T^{\rm H}\bar{u}(\cdot-\ell)$ with $I^{\rm h}$ defined by \eqref{eq:eqI}, then we split it into three parts
\begin{align}
\nonumber
G_2 =& \big \< \big (\delta \FH(T^{\rm H}\bar{u}) - \delta \F^{\rm T}(T^{\rm H}\bar{u})\big)v, v\big \> 
\\[1ex] \nonumber 
=& \sum_{\ell \in \LMM} \big\< \delta \F^{\rm MM}\big(\uh(\cdot-\ell)\big) - \delta T_K \F^{\rm MM}\big(\uh(\cdot-\ell)\big), v \big \> v(\ell) 
\\ \nonumber
& \qquad + \sum_{\ell \in \LMM} \big\< \delta T_K \F^{\rm h}\big(\uh(\cdot-\ell)\big) - \delta \F^{\rm h}\big(\uh(\cdot-\ell)\big), v \big \> v(\ell) 
\\ \nonumber
& + \sum_{\ell \in \LMM} \big\< \delta T_K \FMM\big(\uh(\cdot-\ell)\big) - \delta T_{K}\F^{\rm h}\big(\uh(\cdot-\ell)\big), v \big\> v(\ell)
\\ \nonumber
& + \sum_{\ell \in \LMM} \big \< \delta \F^{\rm h}\big(\uh(\cdot-\ell)\big) - \delta \F^{\rm T}_{\ell}\big(T^{\rm H}\bar{u}\big), v \big\> v(\ell) 
\\[1ex]
=:& G_{21} + G_{22} + G_{23}.
\end{align}
The term $G_{21}$ can be bounded by using the techniques in \cite[Appendix C]{chen17} 
\begin{align}\label{eq:G21_}
|G_{21}| 
\leq C \log\RMM \cdot \RQM^{-K+1/2} \cdot \|Dv\|^2_{\ell^2_{\wf}}.
\end{align}
To estimate $G_{22}$, we have
\begin{eqnarray}\label{eq:G2}
\nonumber
G_{22} 
&=& \sum_{\ell \in \LMM} \big \< \delta \F^{\rm MM}(\p0) - \delta \F^{\rm h}(\p0), v\> v(\ell)
\\ \nonumber
&& + \sum_{j=1}^K \sum_{\ell \in \LMM} \frac{1}{j!}\big \< \delta^{j+1} \F^{\rm MM}(\p0) - \delta^{j+1} \F^{\rm h}(\p0)\big(\uh(\cdot-\ell)\big)^{j}, v\big\> v(\ell) \\[1ex]
&=:& G^{\rm (a)}_{22} + G^{\rm (b)}_{22}.
\end{eqnarray}
We first estimate $G^{\rm (a)}_{22}$ by
\begin{eqnarray}\label{eq:G21}
\nonumber
|G^{\rm (a)}_{22}| 
&\leq& C \ffit_1 \cdot N_{\rm MM}^{1/2} \log\RMM \cdot \|Dv\|^2_{\ell^2_{\wf}}.
\end{eqnarray}
From the above estimate, matching the force constant (i.e. $\ffit_1$ sufficiently small) will lead to the stability of force-mixing schemes \cite{bernstein09, chen17}.  
For $G^{\rm (b)}_{22}$, using Lemma \ref{le:fv} and the decay estimates for the equilibrium state \eqref{eq:ubar-decay}, after a direct calculation, we have
\begin{eqnarray}\label{eq:G1}
|G^{\rm (b)}_{22}| \leq C \sum^{K}_{j=1} \ffit_{j+1} \RQM^{-j} \cdot \log\RMM \cdot \|Dv\|^2_{\ell^2_{\wf}},
\end{eqnarray}
which depends on the measurements $\ffit_j (j\geq2)$ and $\RQM$. 

The term $G_{23}$ can be similarly estimated by \begin{eqnarray}\label{eq:G3}
|G_{23}| \leq C\log\RMM \cdot \big( e^{-\eta \Rbuf} + \RQM^{-K+1/2} \big) \cdot \|Dv\|^2_{\ell^2_{\wf}}.
\end{eqnarray}

Combing from \eqref{eq:stab_hybrid_frc} to \eqref{eq:G3}, for sufficiently large $\RQM$ and sufficiently small $\ffit_1$, we have
\begin{eqnarray}\label{eq:stab_frc}
\big \<\delta \F^{\rm H}(T^{\rm H}\bar{u})v, v\big \>  \geq \frac{c_{\rm F}}{2} \|Dv \|^2_{\ell^2_{\mathcal{N}}}.
\end{eqnarray}

{\it 3. Consistency:} For any $v\in \Admu^{\rm H}$, we have
\begin{eqnarray}\label{eq:proof-c-frc}
\nonumber
&&\<\FH(T^{\rm H}\bar{u}), v\>
\\[1ex]
\nonumber 
&=& \big \<  \F^{\rm T}(T^{\rm H}\bar{u}), v \big\>  + \big \< \FH(T^{\rm H}\bar{u}) - \F^{\rm T}(T^{\rm H}\bar{u}), v\big \> 
\\[1ex]
&=:& P_1 + P_2.
\end{eqnarray}
The term $P_1$ can be estimated by using the results of \cite[Theorem 5.1]{chen17} 
\begin{align}\label{eq:P3}
|P_1|
	\leq C \|Dv \|_{\ell^2_{\wf}}
	\left\{ \begin{array}{ll}
	\RQM^{-(2K+1)d/2}\log\RMM + R^{-d/2}_{\rm MM} + e^{-\kappa_1 \Rbuf} \quad &\text{if}~{\bf (P)},
	\\[1.5ex]
	\RQM^{-K}\log\RMM + \RMM^{-1}\log\RMM + e^{-\kappa_1 \Rbuf} \quad &\text{if}~{\bf (D)}.
	\end{array} \right. 
\end{align}

We then split $P_2$ into three parts
\begin{align}\label{eq:new_P2}
\nonumber
    P_2 =& \big \< \FH(T^{\rm H}\bar{u}) - \F^{\rm T}(T^{\rm H}\bar{u}), v\big \>
    \\[1ex] \nonumber
    =& \sum_{\ell \in \LMM} \Big(\FH\big(\uh(\cdot-\ell)\big) - T_K\FH\big(\uh(\cdot-\ell)\big)\Big) v(\ell)
    \\ \nonumber
    & \qquad + \sum_{\ell \in \LMM} \Big(T_K \F^{\rm h}\big(\uh(\cdot-\ell)\big) - \F^{\rm h}\big(\uh(\cdot-\ell)\big)\Big) v(\ell)
    \\ \nonumber
    &+ \sum_{\ell \in \LMM} \Big(\F^{\rm h}\big(\uh(\cdot-\ell)\big) -  \F_{\ell}^{\rm T}\big(T^{\rm H}\bar{u}\big)\Big) v(\ell)
    \\ \nonumber
    &+ \sum_{\ell \in \LMM} \Big(T_K \FH\big(\uh(\cdot-\ell)\big) - T_K\F^{\rm h}\big(\uh(\cdot-\ell)\big)\Big) v(\ell)
    \\[1ex] 
    =:& P_{21} + P_{22} + P_{23}.
\end{align}
After a direct calculation, the term $P_{21}$ can be estimated by 
\begin{eqnarray}\label{eq:P21_new}
|P_{21}| \leq C\|Dv\|_{\ell^2_{\wf}}
\left\{ \begin{array}{ll}
	\RQM^{-(2K+1)d/2}\log\RMM \quad &\text{if}~{\bf (P)},
	\\[1.5ex]
	\RQM^{-K}\log\RMM \quad &\text{if}~{\bf (D)}.
	\end{array} \right. 
\end{eqnarray}
For $P_{22}$, there exists a constant $\kappa_2>0$ such that
\begin{eqnarray}\label{eq:P22_new}
|P_{22}| \leq C e^{-\kappa_2 \Rbuf}\cdot \|Dv\|_{\ell^2_{\wf}}.
\end{eqnarray}
To estimate $P_{23}$, we further split it into two parts
\begin{align}\label{eq:P1}
\nonumber
P_{23} &= \sum_{\ell \in \LMM} \Big(T_K \FH\big(\uh(\cdot-\ell)\big) - T_K\F^{\rm h}\big(\uh(\cdot-\ell)\big)\Big) v(\ell)
\\[1ex] \nonumber
&= \sum_{\ell \in \LMM} \big( \FH(\p0) - \F^{\rm h}(\p0) \big) \vh(\ell) 
+ \sum_{j=1}^k \frac{1}{j!} \big( \delta^{j} \FH(\p0) - \delta^{j} \F^{\rm h}(\p0)\big)\big(\uh(\cdot-\ell)\big)^{j} \vh(\ell) 
\\[1ex]
&=: P^{\rm (a)}_{23} + P^{\rm (b)}_{23}.
\end{align}
The term $P^{\rm (a)}_{23}$ vanishes since $\FMM(\p0) \equiv \F^{\rm h}(\p0)$. To estimate the second term in \eqref{eq:P1}, we first rewrite it as $P^{\rm (b)}_{23}=:\displaystyle \sum_{j=1}^K \sum_{\ell \in \LMM} \mathscr{F}_j(\ell) v(\ell)$ with 
\vskip-0.3cm
\begin{align}\label{eq:P12}
    |\mathscr{F}_j(\ell)| 
    &\leq \sum_{\pmb{\xi} \in (\Lhom)^j} \left|\frac{ \F^{\rm MM}_{, \pmb{\xi}}(\p0) - \F^{\rm h}_{, \pmb{\xi}}(\p0) }{\prod^{j}_{i=1} \wf(|\xi_i|)} \right| \cdot \prod^{j}_{i=1} \Big(\wf(|\xi_i|) \cdot \uh(\ell + \xi_i) \Big)
    \nonumber \\[1ex] 
    &\leq C \ffit_j \cdot \Big( \sum_{\xi \in \Lhom} \wf(|\xi|) \uh^2(\ell+\xi) \Big)^{j/2} \leq C \ffit_j \cdot  |\ell|^{-\alpha_j},
\end{align}
where $\alpha_j=(d-1)j$ for \asP~or $\alpha_j=j$ for \asD, and the last inequality follows from the decay estimates for the equilibrium state. Hence, using Lemma \ref{le:fv}, the term $P_{23}$ can be bounded by
\begin{align}\label{eq:P23}
|P_{23}| \leq C \|Dv \|_{\ell^2_{\wf}}
	\left\{ \begin{array}{ll}
	\displaystyle
\sum^{K}_{j=1} \ffit_j \RQM^{-(2j-1)d/2} \cdot \log\RMM \quad &\text{if}~{\bf (P)},
	\\[1.5ex]
	\displaystyle
	\sum^{K}_{j=1} \ffit_j \RQM^{-(j-1)} \cdot \log\RMM \quad &\text{if}~{\bf (D)}.
	\end{array} \right. 
\end{align}

Combining from \eqref{eq:proof-c-frc} to \eqref{eq:P23} and let $\kappa:=\max\{\kappa_1, \kappa_2\}$, we have
\begin{align}\label{eq:proof-c-P}
|\<\F^{\rm H}(\THu), v\>| 
\leq
C \|Dv \|_{\ell^2_{\wf}}
	\left\{ \begin{array}{ll}
	\displaystyle
\Big( \sum^{K}_{j=1} \ffit_j \RQM^{-(2j-1)d/2} + \RQM^{-(2K+1)d/2} \Big) \log\RMM \\
\qquad + R^{-d/2}_{\rm MM} + e^{-\kappa \Rbuf} \quad &\text{if}~{\bf (P)},
	\\[1.5ex]
	\displaystyle
	\Big( \sum^{K}_{j=1} \ffit_j \RQM^{-(j-1)} + \RQM^{-K} \Big) \log\RMM +\\
\qquad \RMM^{-1}\log\RMM + e^{-\kappa \Rbuf} \quad &\text{if}~{\bf (D)}.
	\end{array} \right. 
\end{align}
Taking into account stability \eqref{eq:stab_frc} and consistency \eqref{eq:proof-c-P}, we can apply the inverse function theorem \cite[Lemma 2.2]{co2011} with $K=K_{\rm F}$ and use \eqref{proof-TH-diff} to obtain the stated results.
%
\end{proof}

We then give the proof of Theorem \ref{thm:forceM_V}, which is similar to that of Theorem \ref{theo:forceM}. Due to adding the soft matching conditions for virials in Theorem \ref{thm:forceM_V}, the only difference lies in the consistency part.

\begin{proof}[{\bf Proof of Theorem \ref{thm:forceM_V}}]
The difference is the estimate of $P_{21}$ in \eqref{eq:new_P2}. The new estimate of $P_{21}$ now becomes to
\begin{align}\label{eq:P21-v}
\nonumber
P_{21} &= \sum_{\ell \in \LMM} \Big(\FH\big(\uh(\cdot-\ell)\big) - T_K\FH\big(\uh(\cdot-\ell)\big)\Big) v(\ell)
\\ \nonumber
& \qquad + \sum_{\ell \in \LMM} \Big(T_K \F^{\rm h}\big(\uh(\cdot-\ell)\big) - \F^{\rm h}\big(\uh(\cdot-\ell)\big)\Big) v(\ell)
\\[1ex] 
&=: \sum_{\ell \in \LMM} \big( \mathfrak{F}_1\big(\uh(\cdot-\ell)\big) + \mathfrak{F}_2\big(\uh(\cdot-\ell)\big) \big) \cdot v(\ell),
\end{align}
where $\mathfrak{F}_1:= \FH - T_{K+1} \FH + T_{K+1}\F^{\rm h} - \F^{\rm h}$ and $\mathfrak{F}_2:= T_{K+1} \FH - T_{K} \FH + T_{K}\F^{\rm h} -  T_{K+1}\F^{\rm h}$. The first term can be similarly estimated 
\begin{equation}
    \Big| \sum_{\ell \in \LMM} \mathfrak{F}_1 \big( \uh(\cdot-\ell) \big) \cdot v(\ell) \Big| \leq C\|Dv\|_{\ell^2_{\wf}}
\left\{ \begin{array}{ll}
	\RQM^{-(2K+3)d/2}\log\RMM \quad &\text{if}~{\bf (P)},
	\\[1.5ex]
	\RQM^{-(K+1)}\log\RMM \quad &\text{if}~{\bf (D)}.
	\end{array} \right. 
\end{equation}

We expand $\mathfrak{F}_2$ further by using the Taylor's expansion of the strain,
\begin{align}\label{eq:taylor_Du}
\mathfrak{F}_2\big(\uh\big) 
=&~\sum_{(\rho, {\bm \sigma}) \in (\LhomS)^{K+2}} \Big( V^{\rm h}_{, \rho {\bm \sigma}}\big({\bf 0}\big) - V^{\rm MM}_{, \rho {\bm \sigma}}\big({\bf 0}\big) \Big) \Big(\prod^{K+1}_{i=1} \big(\nabla_{\sigma_i}\uh(\ell)+\sum^{\infty}_{t=2}\frac{1}{t!}\nabla^t_{\sigma_i}\uh(\ell)\big) \nonumber \\
&-\prod^{K+1}_{i=1} \big(\nabla_{\sigma_i}u(\ell-\rho)+\sum^{\infty}_{t=2}\frac{1}{t!}\nabla^t_{\sigma_i}u(\ell-\rho)\big)\Big) \nonumber \\[1ex]
=&~\big( \partial^{K+2}_{\rm F} W^{\rm h}_{\rm cb}(\mathsf{I}) - \partial^{K+2}_{\rm F} W^{\rm h}_{\rm cb}(\mathsf{I}) \big) : (\nabla u(\ell))^{K} \nabla^2 u(\ell) + {\rm O}\big(\nabla^2\big(\nabla u(\ell)\big)^{K+1}\big).
\end{align}
Hence, for the second term of \eqref{eq:P21-v}, using the decay estimates for the equilibrium state \eqref{eq:ubar-decay} and Lemma \ref{le:fv}, we have
\begin{equation}\label{eq:P23-v}
    \Big| \sum_{\ell \in \LMM} \mathfrak{F}_2 \big( \uh(\cdot-\ell) \big) \cdot v(\ell) \Big| \leq C \|Dv \|_{\ell^2_{\wf}}
	\left\{ \begin{array}{ll}
\Big( \vfit_{K+1} \RQM^{-(2K+1)d/2} + \RQM^{-(2K+1)d/2-1} \Big) \log\RMM  \quad &\text{if}~{\bf (P)},
	\\[1.5ex]
\Big( \vfit_{K+1} \RQM^{-K} + \RQM^{-K-1} \Big) \log\RMM \quad &\text{if}~{\bf (D)}.
	\end{array} \right. 
\end{equation}
Taking into account the above proofs and the proofs of Theorem \ref{theo:forceM}, we can yield the stated results.
\end{proof}


\section{A semi-empirical QM model: The NRL tight binding}
\label{sec:NRL}
\renewcommand{\theequation}{D.\arabic{equation}}
\setcounter{equation}{0}

In this paper, we use the tight binding model as the reference quantum mechanical model for simplicity of presentation. We note that our numerical scheme is in principle also suitable for general quantum mechanical models. 

The NRL tight binding model is developed by Cohen, Mehl, and Papaconstantopoulos \cite{cohen94}. The energy  levels are determined by the generalised eigenvalue problem
\begin{eqnarray}\label{NRL:diag_Ham}
\mathcal{H}(y)\psi_s = \lambda_s\mathcal{M}(y)\psi_s 
\qquad\text{with}\quad \psi_s^{\rm T}\mathcal{M}(y)\psi_s = 1,
\end{eqnarray}
where $\mathcal{H}$ is the hamiltonian matrix and $\mathcal{M}(y)$ is the overlap matrix. The NRL hamiltonian and overlap matrices are construct both from hopping elements as well as on-site matrix elements as a function of the local environment. For carbon and silicon they are parameterised as follows (for other elements the parameterisation is similar): 

To define the on-site terms, each atom $\ell$ is assigned a pseudo-atomic density 
\begin{eqnarray*}
	\rho_{\ell} := \sum_{k}e^{-\lambda^2 r_{\ell k}}  \fc(r_{\ell k}),
\end{eqnarray*}
where the sum is over all of the atoms $k$ within the cutoff $\Rc$ of atom $\ell$, $\lambda$ is a fitting parameter, $\fc$ is a cutoff function
\begin{eqnarray*}
	\fc(r) = \frac{\theta(\Rc-r)}{1+\exp\big((r-\Rc)/l_c + L_c\big)} ,
\end{eqnarray*}
with $\theta$ the step function, and the parameters $l_c=0.5$, $L_c=5.0$ for most elements.
Although, in principle, the on-site terms should have off-diagonal elements, but this would lead to additional computational challenges that we wished to avoid. The NRL model follows traditional practice and only include the diagonal terms.
Then, the on-site terms for each atomic site $\ell$ are given by
\begin{eqnarray}\label{NRLonsite}
\mathcal{H}(y)_{\ell\ell}^{\upsilon\upsilon}
:= a_{\upsilon} + b_{\upsilon}\rho_{\ell}^{2/3} + c_{\upsilon}\rho_{\ell}^{4/3} + d_{\upsilon}\rho_{\ell}^{2},
\end{eqnarray}
where $\upsilon=s,p$, or $d$ is the index for angular-momentum-dependent atomic orbitals and $(a_\upsilon)$, $(b_\upsilon)$, $(c_\upsilon)$, $(d_\upsilon)$ are fitting parameters.  The on-site elements for the overlap matrix are simply taken to be the identity matrix.

The off-diagonal NRL Hamiltonian entries follow the formalism of  Slater and Koster who showed in \cite{Slater54} that all two-centre (spd) hopping integrals can be constructed from ten independent ``bond integral'' parameters $h_{\upsilon\upsilon'\mu}$, where
\begin{eqnarray*}
	(\upsilon\upsilon'\mu) = ss\sigma,~sp\sigma,~pp\sigma,~pp\pi,
	~sd\sigma,~pd\sigma,~pd\pi,~dd\sigma,~dd\pi,~{\rm and}~dd\delta.
\end{eqnarray*}
The NRL bond integrals are given by
\begin{eqnarray}\label{NRLhopping}
h_{\upsilon\upsilon'\mu}(r) 
:= \big(e_{\upsilon\upsilon'\mu} + f_{\upsilon\upsilon'\mu}r + g_{\upsilon\upsilon'\mu} r^2 \big) e^{-h_{\upsilon\upsilon'\mu}r} \fc(r)
\end{eqnarray}
with fitting parameters $e_{\upsilon\upsilon'\mu}, f_{\upsilon\upsilon'\mu},  g_{\upsilon\upsilon'\mu}, h_{\upsilon\upsilon'\mu}$. The matrix elements $\mathcal{H}(y)_{\ell k}^{\upsilon\upsilon'}$ are constructed from the $h_{\upsilon\upsilon'\mu}(r)$ by a standard procedure~\cite{Slater54}.

The analogous bond integral parameterisation of the overlap matrix 
is given by  
\begin{eqnarray}\label{NRLhopping-M}
m_{\upsilon\upsilon'\mu}(r) 
:= \big(\delta_{\upsilon\upsilon'} + p_{\upsilon\upsilon'\mu} r + q_{\upsilon\upsilon'\mu} r^2 + r_{\upsilon\upsilon'\mu} r^3 \big) e^{-s_{\upsilon\upsilon'\mu} r} \fc(r)
\end{eqnarray}
with the fitting parameters $(p_{\upsilon\upsilon'\mu}), (q_{\upsilon\upsilon'\mu}), (r_{\upsilon\upsilon'\mu}), (s_{\upsilon\upsilon'\mu})$ and $\delta_{\upsilon\upsilon'}$ the Kronecker delta function. 

The fitting parameters in the foregoing expressions are determined by fitting to some high-symmetry first-principle calculations: In the NRL method, a database of eigenvalues (band structures) and total energies were constructed for several crystal structures at  several volumes. Then the parameters are chosen such that the eigenvalues and energies in the database are reproduced.
For practical simulations, the parameters for different elements can be found in \cite{papaconstantopoulos15}.


\section{The Atomic Cluster Expansion}
\label{sec:ACE}
\renewcommand{\theequation}{E.\arabic{equation}}
\setcounter{equation}{0}


Following \cite{bachmayr19}, we briefly introduce the construction of the ACE potential. 
Given $\mathcal{N} \in \N$, we first write the ACE site potential in the form of an {\it atomic body-order expansion}, $\displaystyle V^{\rm ACE} \big(\{\pmb{g_j}\}\big) = \sum_{N=0}^{\mathcal{N}} \frac{1}{N!}\sum_{j_1 \neq \cdots \neq j_N} V_N(\pmb{g}_{j_1}, \cdots, \pmb{g}_{j_N})$, where the $N$-body potential $V_N : \R^{dN} \rightarrow \R$ can be approximated by using a tensor product basis \cite[Proposition 1]{bachmayr19},
\begin{align*}
\phi_{\pmb{k\ell m}}\big(\{\pmb{g}_j\}_{j=1}^N\big) :=\prod_{j=1}^N\phi_{k_j\ell_j m_j}(\pmb{g}_j) 
\quad & {\rm with}\quad 
\phi_{k\ell m}(\pmb{r}):=P_k(r)Y^{m}_{\ell}(\hat{r}) , ~~\pmb{r}\in\R^d,~r=|\pmb{r}|,~\hat{r}=\pmb{r}/r ,
\end{align*}
where $P_k,~k=0,1,2,\cdots$ are radial basis functions (for example, Jacobic polynomials), and $Y_{\ell}^m,~\ell = 0,1,2,\cdots,~m=-\ell,\cdots,\ell$ are the complex spherical harmonics.
The basis functions are further symmetrised to a permutation invariant form,
\begin{eqnarray*}
\tilde{\phi}_N = \sum_{(\pmb{k,\ell,m})~{\rm ordered}} \sum_{\sigma\in S_N} \phi_{\pmb{k\ell m}} \circ \sigma ,
\end{eqnarray*}  
where $S_N$ is the collection of all permutations, and by $\sum_{(\pmb{k,\ell,m})~{\rm ordered}}$ we mean that the sum is over all lexicographically ordered tuples $\big( (k_j, \ell_j, m_j) \big)_{j=1}^N$.
The next step is to incorporate the invariance under point reflections and rotations
\begin{eqnarray*}
\mathcal{B}_{\pmb{k\ell} i} = \sum_{\pmb{m}\in\mathcal{M}_{\pmb{\ell}}}\mathcal{U}_{\pmb{m}i}^{\pmb{k\ell}} \sum_{\sigma\in S_N} \phi_{\pmb{k\ell m}} \circ \sigma 
\quad{\rm with}\quad \mathcal{M}_{\pmb{\ell}} = \big\{\pmb{\mu}\in\Z^N ~|~ -\ell_{\alpha}\leq\mu_{\alpha}\leq\ell_{\alpha} \big\} ,
\end{eqnarray*}
where the coefficients $\mathcal{U}_{mi}^{k\ell}$ are given in \cite[Lemma 2 and Eq. (3.12)]{bachmayr19}.
It was shown in \cite{bachmayr19} that the basis defined above is explicit but computational inefficient. The so-called ``density trick" technique used in \cite{Bart10, Drautz19, Shapeev16} can transform this basis into one that is computational efficient. The alternative basis is 
\begin{eqnarray*}
B_{\pmb{k\ell} i} = \sum_{\pmb{m}\in\mathcal{M}_{\pmb{\ell}}}\mathcal{U}_{\pmb{m}i}^{\pmb{k\ell}} A_{\pmb{nlm}}
\quad{\rm with~the~correlations}\quad A_{\pmb{nlm}}:= \prod_{\alpha=1}^{N} \sum_{j=1}^{J} \phi_{n_\alpha l_\alpha m_{\alpha}}(\pmb{g}_j),
\end{eqnarray*}
which avoids both the $N!$ cost for symmetrising the basis as well as the $C_J^N$ cost of summation over all order $N$ clusters within an atomic neighbourhood.
The resulting basis set is then defined by
\begin{eqnarray}
\pmb{B}_N := \big\{ B_{\pmb{k\ell} i}  ~|~ (\pmb{k},\pmb{\ell})\in\N^{2N}~{\rm ordered}, ~\sum_{\alpha} \ell_{\alpha}~{\rm even},~ i=1,\cdots,\pmb{n}_{\pmb{k\ell}}\big\},
\end{eqnarray}
where $\pmb{n}_{\pmb{k\ell}}$ is the rank of body-orders (see \cite[Proposition 7 and Eq. (3.12)]{bachmayr19}).

Once the finite symmetric polynomial basis set $\pmb{B} \subset \bigcup^{\mathcal{N}}_{N=1} \pmb{B}_N$ is constructed, the ACE site potential can be expressed as
\begin{align}\label{ships:energy}
V^{\rm ACE}(\pmb{g}; \{c_B\}_{B\in\pmb{B}}) = \sum_{B\in\pmb{B}} c_B B(\pmb{g}) 
\end{align}
with the coefficients $c_{B}$. The corresponding force of this potential is denoted by $\F^{\rm ACE}$.

The family of potentials are systematically improvable (see \cite[\S 6.2]{bachmayr19}): by increasing the body-order, cutoff radius and polynomial degree they are in principle capable of representing an arbitrary many-body potential energy surface to within arbitrary accuracy.


\section{Numerical supplements}
\label{sec:ns}
\renewcommand{\theequation}{F.\arabic{equation}}
\setcounter{equation}{0}
\renewcommand{\thetable}{F.\arabic{table}}
\setcounter{table}{0}

\begin{table}
	\begin{center}
    {\bf Weights for energy-mixing}\\
	\vskip0.2cm
	\begin{tabular}{|c|c|c|c|c|c|}
	\hline
		\diagbox{Schemes}{Weights} & {~~~$W^{\rm E}_0$~~~} & {~~~$W^{\rm E}_1$~~~} & {~~~$W^{\rm E}_2$~~~} & {~~~$W^{\rm E}_3$~~~} & {~~~$W^{\rm V}_2$~~~}\\[0.5mm]
		\hline
		$K_{\rm E}=2$  &  1  & 10 & 500 & - & -   \\[0.5mm]
		\hline
		$K_{\rm E}=3$  &  1  & 10 & 1000 & 100 & -
	    \\[0.5mm]
	    \hline
		$K_{\rm E}=2 ~~\&$  virial    &  1  & 10 & 1000 & - & 100
		\\[0.5mm]
		\hline
	\end{tabular}
	\\[6mm]
    {\bf Weights for force-mixing}\\
	\vskip0.2cm
	\begin{tabular}{|c|c|c|c|c|c|}
	\hline
		\diagbox{Schemes}{Weights}& {~~~$W^{\rm F}_0$~~~} & {~~~$W^{\rm F}_1$~~~} & {~~~$W^{\rm F}_2$~~~} & {~~~$W^{\rm V}_2$~~~} & {~~~$W^{\rm V}_3$~~~}\\[0.5mm]
		\hline
		$K_{\rm F}=1$  &  1  & 100 & - & - & -   \\[0.5mm]
		\hline
		$K_{\rm F}=2$  &  1  & 1000 & 100 & - & -
	    \\[0.5mm]
	    \hline
		$K_{\rm F}=1~~\&$ virial    &  1  & 1000 & - & 100 & -
		\\[0.5mm]
	   \hline
	     $K_{\rm F}=2~~\&$ virial    &  1  & 1000 & 200 & - & 500
		\\[0.5mm]
		\hline
	\end{tabular}
	\end{center}
	\vskip 1cm
	\caption{The practical choice of weights in the loss functional for energy-mixing (varying $K_{\rm E}$) and force-mixing (varying $K_{\rm F}$) coupling schemes.
	}
	\label{table-weights}
\end{table}

\small
\bibliographystyle{plain}
\bibliography{bib}

\end{document}